\documentclass[a4paper,draft,reqno]{amsart}
\usepackage{mathptmx}
\usepackage{amssymb}
\usepackage{amsmath}
\usepackage{amsfonts}
\usepackage{amscd}
\usepackage{enumerate}
\usepackage{epsfig}
\usepackage{amscd}
\newcommand{\lattice}{\mathcal{S}}

\newcommand{\field}{\mathcal{F}}

\newcommand{\pspace}{\Omega}
\newcommand{\gambles}{\mathcal{L}}
\newcommand{\linlat}{\mathcal{K}}

\newcommand{\solp}{\mathcal{M}}

\newcommand{\meet}{\wedge}
\newcommand{\join}{\vee}
\newcommand{\pr}{P}
\newcommand{\lpr}{{\underline{P}}}
\newcommand{\upr}{{\overline{P}}}
\newcommand{\lnex}{{\underline{E}}}

\newcommand{\nex}{E}
\newcommand{\apr}{Q}
\newcommand{\alpr}{{\underline{Q}}}

\newcommand{\fnl}{\Gamma} 
\newcommand{\lowfnl}{{\underline{\fnl}}} 
\newcommand{\uppfnl}{{\overline{\fnl}}} 
\newcommand{\linfnl}{\fnl} 
\newcommand{\linfnls}{{\gambles^*_+}} 
\newcommand{\afnl}{\Psi} 
\newcommand{\alowfnl}{{\underline{\afnl}}} 
\newcommand{\alinfnl}{\afnl} 
\newcommand{\aafnl}{\Xi} 
\newcommand{\aalowfnl}{{\underline{\aafnl}}} 
\newcommand{\norm}[1]{{\lVert#1\rVert}}
\newcommand{\card}[1]{{\lvert#1\rvert}}
\newcommand{\set}[2]{\left\{#1\colon#2\right\}}
\newcommand{\ddf}[2]{{G^{#2}_{#1}}}
\newcommand{\dif}{{\,\mathrm{d}}}
\DeclareMathOperator{\cl}{cl}
\DeclareMathOperator{\linspan}{span}

\newcommand{\SetN}{{\mathbb{N}}}
\newcommand{\SetNinfty}{\SetN^*}
\newcommand{\SetNzero}{\SetN_0}
\newcommand{\SetR}{\mathbb{R}}
\DeclareMathOperator{\Domain}{dom}

\newcommand{\Abs}[1]{\left|{#1}\right|}

\theoremstyle{plain}
\newtheorem{theorem}{Theorem}
\newtheorem{proposition}[theorem]{Proposition}
\newtheorem{corollary}[theorem]{Corollary}
\newtheorem{lemma}[theorem]{Lemma}
\theoremstyle{remark}
\newtheorem{example}{Example}
\newtheorem{counter}{Counterexample}

\theoremstyle{definition}
\newtheorem{definition}{Definition}
\begin{document}
\title{$n$-Monotone exact functionals}
\author{Gert de Cooman}
\address{Ghent University, SYSTeMS Research Group, Technologiepark -- Zwijnaarde 914, 9052
  Zwijnaarde, Belgium. E-mail: \texttt{gert.decooman@ugent.be}.}
\author{Matthias C.~M.~Troffaes}
\address{Durham University, Dep. of Mathematical Sciences,  Science Laboratories, South Road, Durham, DH1 3LE, England. E-mail: \texttt{matthias.troffaes@gmail.com}.}
\author{Enrique Miranda}
\address{Rey Juan Carlos University, Dep. of Statistics and O. R., C-Tulip\'an,
s/n, 28933, M\'ostoles, Spain. E-mail: \texttt{enrique.miranda@urjc.es}.}
\begin{abstract}
  We study $n$-monotone functionals, which constitute a generalisation of $n$-monotone set functions.   We investigate their relation to the concepts of exactness and natural extension, which generalise   the notions of coherence and natural extension in the behavioural theory of imprecise   probabilities. We improve upon a number of results in the literature, and prove among other things a   representation result for exact $n$-monotone functionals in terms of Choquet integrals.
\end{abstract}
\date{June 2006}
\maketitle

\section{Introduction}

Exact functionals are real-valued functionals that are monotone, super-additive, positively homogenous, and translation invariant (or constant additive). They were introduced and studied by Maa\ss\ \cite{maass2002,maass2003} in an attempt to unify and generalise a number of notions in the literature, such as coherent lower previsions (Walley \cite{walley1991}), exact cooperative games (Schmeidler \cite{schmeidler1972}) and coherent risk measures (Artzner \textit{et al.}  \cite{artzner1999}, Delbaen \cite{delbaen2002}).

Coherent lower previsions, mainly due to Walley \cite{walley1991}, are among the most interesting uncertainty models in what has been called the theory of imprecise probabilities; this is the theory which extends the Bayesian theory of probability by allowing for indecision. Coherent lower previsions can be viewed as lower expectations with respect to closed convex sets of probability measures (also called credal sets; see Levi \cite{levi1980a}), and they provide a unifying framework for studying many other uncertainty models, such as probability charges (Bhaskara Rao and Bhaskara Rao \cite{bhaskara1983}), $2$- and $n$-monotone set functions (Choquet \cite{choquet1953}), possibility measures (\cite{cooman2001,cooman1998,cooman1998a,dubois1988}), and p-boxes (Ferson \textit{et al.} \cite{2003:ferson::pboxes}). They have also been linked to various theories of integration, such as Choquet integration (Walley \cite[p.~53]{walley1981}) and Lebesgue integration (Walley \cite[p.~132]{walley1991}). Exact functionals are essentially coherent lower previsions multiplied by a non-negative constant (see Theorem~\ref{theo:exactness-reduces-to-coherence} further on). On the other hand, the coherent risk measures introduced by Artzner \textit{et   al.}  (\cite{artzner1999,delbaen2002}), which have become quite important in finance theory, are just the negatives of exact functionals.

Here, we study the properties of a special subclass of exact functionals, namely those that are \emph{$n$-monotone}, for $n\geq1$. We start out from Choquet's \cite{choquet1953} original and very general definition of $n$-mono\-toni\-city for functions defined on arbitrary lattices, and we pave the way towards a representation theorem for $n$-monotone exact functionals in terms of the Choquet integral.

The paper is structured as follows.  Section~\ref{sec:lowprevs:shortintro} highlights the most important aspects of the theory of coherent lower previsions needed in the rest of the paper, and Section~\ref{sec:exactfunctionals} explains their generalisation to exact functionals. Section~\ref{sec:definition} is concerned with the precise definition of $n$-monotonicity for exact functionals. In Section~\ref{sec:lower-probabilities}, we establish many interesting properties, and generalise a number of results from the literature for $n$-monotone set functions on fields of events. In Section~\ref{sec:representation}, we relate $n$-monotone exact functionals to comonotone additive functionals and Choquet integrals. We conclude in Section~\ref{sec:conclusions} with some additional comments and remarks.

\section{Coherent lower previsions: a short introduction}
\label{sec:lowprevs:shortintro}
In this section, we introduce a few basic notions about coherent lower previsions. We refer to Walley \cite{walley1991} for a more in-depth discussion.

Consider a non-empty set $\pspace$. A \emph{gamble} $f$ on $\pspace$ is a bounded real-valued mapping on $\pspace$. The set of all gambles on $\pspace$ is denoted by $\gambles$. It is a real linear space under the point-wise addition of gambles, and the point-wise scalar multiplication of gambles with real numbers. Given a real number $\mu$, we also use $\mu$ to denote the gamble that takes the constant value $\mu$.

A special class of gambles are the ones that only take values in $\{0,1\}$: let $A$ be any subset of $\pspace$, also called an \emph{event}, then the gamble $I_A$, defined by $I_A(\omega):=1$ if $\omega\in A$ and $I_A(\omega):=0$ otherwise, is called the \emph{indicator} of $A$. This establishes a correspondence between events and $\{0,1\}$-valued gambles.  Often, for an event $A$, we also denote $I_A$ by $A$.

A \emph{lower prevision} $\lpr$ is a real-valued map (a functional) defined on some subset of $\gambles$, called its \emph{domain} and denoted by $\Domain\lpr$. For any gamble $f$ in $\Domain\lpr$, $\lpr(f)$ is called the lower prevision of $f$.  If the domain of $\lpr$ contains only (indicators of) events $A$, then we also call $\lpr$ a \emph{lower probability}, and we write $\lpr(I_A)$ also as $\lpr(A)$, the lower probability of $A$.

Given a lower prevision $\lpr$, its \emph{conjugate upper prevision} $\upr$ is defined on the set of gambles $\Domain\upr=-\Domain\lpr:=\set{-f}{f\in\Domain\lpr}$ by $\upr(f):=-\lpr(-f)$ for every $-f$ in the domain of $\lpr$. This conjugacy relationship shows that we can restrict our attention to the study of lower previsions only. If the domain of $\upr$ contains indicators only, then we also call $\upr$ an \emph{upper probability}.

Recall that a linear space of gambles is a subset of $\gambles$ that is closed under point-wise addition of gambles and scalar multiplication of gambles with real numbers. Then a lower prevision $\lpr$ whose domain is a linear space is called \emph{coherent} if the following three properties are satisfied for all $f$, $g$ in $\Domain\lpr$ and all non-negative real $\lambda$:
\begin{enumerate}[(C1)]
\item $\lpr(f)\geq\inf f$ (accepting sure gains);
\item $\lpr(\lambda f)=\lambda\lpr(f)$ (positive homogeneity);
\item $\lpr(f+g)\geq\lpr(f)+\lpr(g)$ (super-additivity).
\end{enumerate}
It can be shown that a coherent lower prevision on a linear space can always be extended to a coherent lower prevision on all gambles.

A lower prevision $\lpr$ with a general domain (not necessarily a linear space) is then called \emph{coherent} if it can be extended to a coherent lower prevision on all gambles.  This is the case if and only if $\sup\left[\sum_{i=1}^nf_i-mf_0\right]\geq\sum_{i=1}^n\lpr(f_i)-m\lpr(f_0)$ for any natural numbers $n\geq0$ and $m\geq0$, and $f_0$, $f_1$, \dots, $f_n$ in the domain of $\lpr$.

There are a number of common consequences of coherence that we shall use further on.  Consider a coherent lower prevision $\lpr$, let $f$ and $g$ be elements in $\Domain\lpr$, and let $\mu$ and $\lambda$ be real numbers, with $\lambda\geq0$. Then whenever the relevant gambles belong to $\Domain\lpr$, we have that $\lpr(f+g)\geq\lpr(f)+\lpr(g)$, $\lpr(\lambda f)=\lambda\lpr(f)$, $\lpr(\mu)=\mu$ and $\lpr(f+\mu)=\lpr(f)+\mu$. Moreover $\inf f\leq\lpr(f)\leq\upr(f)\leq\sup f$ and consequently $0\leq\lpr(\Abs{f})\leq\upr(\Abs{f})\leq\sup\Abs{f}$. Also, $\lpr$ is monotone: if $f\leq g$ then $\lpr(f)\leq\lpr(g)$. Finally, both $\Abs{\lpr(f)-\lpr(g)}\leq\upr(\Abs{f-g})$ and $\Abs{\upr(f)-\upr(g)}\leq\upr(\Abs{f-g})$. As an immediate consequence of these properties, we see that if a sequence $f_n$ of gambles converges uniformly to a gamble $f$, \textit{i.e.}, $\sup\Abs{f_n-f}\to0$, then also $\lpr(f_n)\to\lpr(f)$ and $\upr(f_n)\to\upr(f)$, so any coherent lower or upper prevision is continuous with respect to the supremum norm.

A lower prevision $\alpr$ is said to \emph{dominate} a lower prevision $\lpr$, if $\Domain\alpr\supseteq\Domain\lpr$ and $\alpr(f)\ge\lpr(f)$ for any $f$ in $\Domain\lpr$. We say that a lower prevision $\lpr$ \emph{avoids sure loss} if it is dominated by some coherent lower prevision on $\gambles$. This is the case if and only if $\sup\left[\sum_{i=1}^nf_i\right]\ge\sum_{i=1}^n\lpr(f_i)$ for any natural number $n\geq1$ and any $f_1$, \dots, $f_n$ in $\Domain\lpr$.

One can easily show that a lower prevision avoids sure loss if and only if there is a point-wise smallest coherent lower prevision $\lnex_\lpr$ on $\gambles$ that dominates $\lpr$, namely, the lower envelope of all the coherent lower previsions on $\gambles$ that dominate $\lpr$ on $\Domain\lpr$.  $\lnex_\lpr$ is then called the \emph{natural extension} of $\lpr$. It is also given by (Walley \cite[Lemma~3.1.3(b)]{walley1991})
\begin{multline}\label{eq:nextension}
  \lnex_\lpr(f)
  =\sup\Bigg\{\sum_{k=1}^n\lambda_k\lpr(f_k)+\lambda\colon \\*
  n\geq1,\lambda_k\in\SetR_+,\lambda\in\SetR,f_k\in\Domain\lpr,   \sum_{k=1}^n\lambda_kf_k+\lambda\leq f\Bigg\}
\end{multline}
for all $f \in \gambles$, where $\SetR_+$ is the set of non-negative real numbers.

A \emph{linear prevision} $\pr$ is a real-valued functional defined on a set of gambles $\Domain\pr$, that satisfies $\sup[\sum_{i=1}^nf_i-\sum_{j=1}^mg_j] \geq\sum_{i=1}^n\pr(f_i)-\sum_{j=1}^m\pr(g_j)$ for any natural numbers $n$ and $m$, and $f_1$, \dots, $f_n$, $g_1$, \dots, $g_m$ in the domain of $\lpr$. Note that a linear prevision $\pr$ is coherent, both when interpreted as a lower, and as an upper prevision; the former means that $\pr$ is a coherent lower prevision on $\Domain\pr$, the latter that $-\pr(-\cdot)$ is a coherent lower prevision on $-\Domain\pr$. For any linear prevision $\pr$, it holds that $\pr(f)=-\pr(-f)$ whenever $f$ and $-f$ belong to the domain of $\pr$. A lower prevision $\lpr$ whose domain is negation invariant (\textit{i.e.}, $-\Domain\lpr=\Domain\lpr$), is a linear prevision if and only if it is coherent and \emph{self-conjugate}; self-conjugacy means that $\lpr(-f)=-\lpr(f)$ for all $f$ in $\Domain\lpr$.  A linear prevision $\pr$ on $\gambles$ is easily seen to be a non-negative, normed [$\pr(1)=1$], real-valued, linear functional on $\gambles$.  The restriction of such a linear prevision on $\gambles$ to (indicators of) events is a \emph{probability charge} (or finitely additive probability measure) on $\wp(\pspace)$, the class of all subsets of $\pspace$.

Let us denote the set of linear previsions on $\gambles$ that dominate $\lpr$ by $\solp(\lpr)$. The following statements are equivalent: (i) $\lpr$ avoids sure loss, (ii) the natural extension of $\lpr$ exists; and (iii) $\solp(\lpr)$ is non-empty.  The following statements are equivalent as well: (i) $\lpr$ is coherent; (ii) $\lpr$ coincides with its natural extension $\lnex_\lpr$ on $\Domain\lpr$; and (iii) $\lpr$ coincides with the lower envelope of $\solp(\lpr)$ on $\Domain\lpr$.  The last statement follows from the important fact that $\lnex_\lpr$ 
is equal to the lower envelope of $\solp(\lpr)$:
\begin{equation*}
  \lnex_\lpr(f)=\min_{\apr\in\solp(\lpr)}\apr(f),
\end{equation*}
for any gamble $f$ in $\gambles$. Often, this expression provides a convenient way of calculating the natural extension of a lower prevision that avoids sure loss. Finally it holds that $\solp(\lpr)=\solp(\lnex_\lpr)$. This result can be used to prove the following ``transitivity'' property for natural extension: if we denote by $\alpr$ the restriction of the natural extension $\lnex_\lpr$ of a lower prevision (that avoids sure loss) to some set of gambles $\linlat\supseteq\Domain\lpr$, then $\solp(\lpr)=\solp(\alpr)=\solp(\lnex_\lpr)$, and consequently $\lnex_\alpr$ coincides with $\lnex_\lpr$.

\section{Exact functionals}
\label{sec:exactfunctionals}

\subsection{Notation and Definitions}

In what follows, we use the term \emph{functional} to refer to a real-valued map defined on some subset of $\gambles$.  If $\lowfnl$ denotes a functional, then $\uppfnl$ denotes its \emph{conjugate}, defined by
\begin{equation*}
  \uppfnl(f):=-\lowfnl(-f),
\end{equation*}
for any gamble $f$ in $-\Domain\lowfnl:=\{-f\colon f\in\Domain\lowfnl\}$. So, $\Domain\uppfnl=-\Domain\lowfnl$.

Maa\ss\ \cite{maass2003} has extended the notion of coherence for lower previsions to that of \emph{ exactness} for functionals: a functional $\lowfnl$ on $\gambles$ is called \emph{exact} whenever for any gambles $f$ and $g$ on $\pspace$, any non-negative real number $\lambda$, and any real number $\mu$, it holds that
\begin{enumerate}[(E1)]
\item if $f\geq g$ then $\lowfnl(f)\geq\lowfnl(g)$ (monotonicity);
\item $\lowfnl(\lambda f)=\lambda\lowfnl(f)$ (positive homogeneity);
\item $\lowfnl(f+g)\geq\lowfnl(f)+\lowfnl(g)$ (super-additivity);
\item $\lowfnl(f+\mu)=\lowfnl(f)+\lowfnl(\mu)$ (constant additivity).
\end{enumerate}
A functional defined on an arbitrary subset of $\gambles$ is called \emph{exact} if it can be extended to an exact functional on all of $\gambles$.

The conjugates of exact functionals generalise coherent upper previsions, and are sub-additive rather than super-additive.

An exact functional $\linfnl$, defined on an arbitrary subset of $\gambles$, is called \emph{linear} if it can be extended to an exact functional $\alinfnl$ on $\gambles$ which is at the same time a linear functional, \textit{i.e.}, which also satisfies $\alinfnl(f)+\alinfnl(g)=\alinfnl(f+g)$ for any $f$ and $g$ in $\gambles$. The linear exact functionals on $\gambles$ are precisely the positive linear functionals on $\gambles$. We denote the set of all linear exact functionals on $\gambles$ by $\linfnls$.

Let us give a simple example of a non-exact positive linear functional. Here and elsewhere in this paper the set of natural numbers without zero is denoted by $\SetN$. By $\SetNinfty$ we denote $\SetN\cup\{\infty\}$ and by $\SetNzero$ the set $\SetN\cup\{0\}$. Consider the linear space $\mathcal{K}$, defined by
\begin{equation*}
  \mathcal{K}
  :=\left\{
    \sum_{i=1}^n y_i I_{[a_i,b_i]}\colon
    n\in\SetN,y_1,a_1,b_1,\dots,y_n,a_n,b_n\in\SetR,\,a_1<b_1,\dots,a_n<b_n
  \right\}.
\end{equation*}
Define $\linfnl$ on $\mathcal{K}$ as the Lebesgue-integral on $\mathcal{K}$:
\begin{equation*}
  \linfnl\left(\sum_{i=1}^n y_i I_{[a_i,b_i]}\right):=\sum_{i=1}^n y_i(b_i-a_i).
\end{equation*}
Clearly, $\linfnl$ is a real-valued, and it is a positive linear functional.  But it is not exact, simply because it is not continuous with respect to the supremum norm, and such continuity is a property that all exact functionals have, as we shall see at the end of Section~\ref{sec:exactness-and-coherence}: even though the sequence of gambles $f_n:=\frac{1}{n}I_{[0,n]}$ converges uniformly to the zero gamble, $\linfnl(f_n)=\frac{1}{n}(n-0)=1$ does not converge to zero. This also proves that $\linfnl$ has no exact extension to the set of all gambles on $\SetR$.

It can be proven that a positive linear functional $\linfnl$ on a linear lattice $\mathcal{K}$ is exact if and only if $\linfnl$ is continuous with respect to the supremum norm. However, the equivalence does not necessarily hold if the domain $\mathcal{K}$ is not a linear lattice of gambles.

\subsection{The relation between exactness and coherence}
\label{sec:exactness-and-coherence}
Consider an exact functional $\lowfnl$, then clearly for any $\lambda\geq0$ the functional $\lambda\lowfnl$ is exact as well.  Moreover, if a functional $\lowfnl$ is exact, and both $\mu\in\SetR$ and $1$ belong to its domain $\Domain\lowfnl$, then it follows easily that
\begin{enumerate}
\item[(E5)] $\lowfnl(\mu)=\mu\lowfnl(1)$.
\end{enumerate}
Therefore, a coherent lower prevision $\lpr$, whose domain contains at least the constant gamble $1$, is an exact functional which additionally satisfies $\lpr(1)=1$.  We shall see further on in Theorem~\ref{theo:exactness-reduces-to-coherence} that exact functionals are essentially coherent lower previsions, but without the normalisation constraint $\lpr(1)=1$.

To obtain this result, we use the following norm defined on functionals, introduced by Maa\ss\ \cite[Eq.~(1.2), p.~4]{maass2003}:
\begin{equation*}
  \norm{\lowfnl}
  :=\inf\set{c\in\SetR_+}%
  {f\geq\sum_{k=1}^n\lambda_kf_k+\lambda\Rightarrow
    \lowfnl(f)\geq\sum_{k=1}^n\lambda_k\lowfnl(f_k)+\lambda c},
\end{equation*}
where the condition must hold for all $n$ in $\SetN$, $\lambda_1$, \dots, $\lambda_n$ in $\SetR_+$, $\lambda$ in $\SetR$, and gambles $f$, $f_1$, \dots, $f_n$ in $\Domain\lowfnl$. It holds that $\norm{\lowfnl}\ge 0$ and $\norm{\lowfnl}=0$ implies $\lowfnl=0$, $\norm{\lambda\lowfnl}=\lambda\norm{\lowfnl}$, and $\norm{\lowfnl+\alowfnl}\le\norm{\lowfnl}+\norm{\alowfnl}$, for any functionals $\lowfnl$ and $\alowfnl$ defined on the same domain, and any non-negative real $\lambda$ (see Maa\ss\ \cite[Prop.~1.2.3(a)--(c)]{maass2003}); this motivates our calling $\norm{\lowfnl}$ the norm of $\lowfnl$.

Maa{\ss} \cite[Prop.~1.2.4]{maass2003} has proven that if $\lowfnl$ is an exact functional such that $1\in\Domain\lowfnl$, then $\norm{\lowfnl}=\lowfnl(1)$; this yields a convenient expression for the norm.  He has also proven the following theorem, which shows that exactness of a functional $\lowfnl$ is completely determined by its norm $\norm{\lowfnl}$, and which provides us with a constructive way to obtain an exact extension of $\lowfnl$ to the set $\gambles$ of all gambles on $\pspace$, similar to natural extension for lower previsions.

\begin{theorem}[\protect{Maa\ss\ \cite[Thm.~1.2.5]{maass2003}}]
  \label{theo:maass-results}
  Any functional\/ $\lowfnl$ is exact if and only if\/ $\norm{\lowfnl}<+\infty$. Moreover,   if\/ $\lowfnl$ is exact then the functional\/ $\lnex_\lowfnl$ on $\gambles$ defined for all   gambles $f$ on $\pspace$ by
  \begin{equation*}
    \lnex_\lowfnl(f)
    =\sup\set{\sum_{k=1}^n\lambda_k\lowfnl(f_k)+\lambda\norm{\lowfnl}}
    {\sum_{k=1}^n\lambda_kf_k+\lambda\leq f},
  \end{equation*}
  where the supremum runs over all $n$ in $\SetN$, $\lambda_1$, \dots, $\lambda_n$ in   $\SetR_+$, $\lambda$ in $\SetR$, and gambles $f_1$, \dots, $f_n$ in $\Domain\lowfnl$, is an   exact extension of\/ $\lowfnl$ with $\norm{\lnex_\lowfnl}=\norm{\lowfnl}=\lnex_\lowfnl(1)$.
\end{theorem}

An exact functional $\lowfnl$ has by definition exact extensions to all of $\gambles$. We now see that it also has at least one exact extension $\lnex_\lowfnl$ whose norm is equal to $\norm{\lowfnl}$.  We can associate with $\lowfnl$ its set of dominating positive linear functionals on $\gambles$ with the same norm:
\begin{equation*}
  \solp(\lowfnl)
  :=\set{\alinfnl\in\linfnls}
  {\text{$\alinfnl\geq\lowfnl$ and $\norm{\alinfnl}=\norm{\lowfnl}$}},
\end{equation*}
where $\alinfnl\geq\lowfnl$ means that $\alinfnl(f)\geq\lowfnl(f)$ for every gamble $f$ on $\Domain\lowfnl$. Then $\lnex_\lowfnl$ is the lower envelope of $\solp(\lowfnl)$ and moreover $\solp(\lowfnl)=\solp(\lnex_\lowfnl)$. These results follow at once from Theorem~\ref{theo:exactness-reduces-to-coherence} below, and the corresponding results mentioned in the previous section for coherent lower previsions. An alternative proof can be found in Maa{\ss} \cite[Prop.~1.2.7]{maass2003}. The exact functional $\lnex_\lowfnl$ is called the \emph{natural extension} of the exact functional $\lowfnl$. Just like its counterpart for coherent lower previsions, the natural extension of exact functionals is ``transitive'' (see the discussion at the end of Section~\ref{sec:lowprevs:shortintro}).

We now prove a theorem that uncovers the relationship between coherent lower previsions, exact functionals, and their natural extensions.

\begin{theorem}\label{theo:exactness-reduces-to-coherence}
  Let\/ $\lowfnl$ be a functional defined on a subset of $\gambles$. The following holds.
  \begin{enumerate}[(i)]
  \item\label{theo:exactness-reduces-to-coherence:simple} If\/ $\lowfnl$ is exact, then there     is a coherent lower prevision\/ $\lpr$ defined on\/ $\Domain\lowfnl$ such that\/     $\lowfnl=\norm{\lowfnl}\lpr$, and moreover\/ $\lnex_\lowfnl=\norm{\lowfnl}\lnex_\lpr$.
  \item\label{theo:exactness-reduces-to-coherence:complex} $\lowfnl$ is exact if and only if     there is a coherent lower prevision $\lpr$ defined on $\Domain\lowfnl$, and a non-negative     real number $\lambda$, such that\/ $\lowfnl=\lambda\lpr$.  In that case,     $\lambda\lnex_\lpr$ is an exact extension of\/ $\lowfnl$ with norm $\lambda$.

    If, additionally, $1$ belongs to the domain of\/ $\lowfnl$, then $\lambda$ is uniquely     given by $\lowfnl(1)$, and hence, $\lnex_\lowfnl=\lowfnl(1)\lnex_\lpr$; and if also     $\lowfnl$ is non-zero for at least one gamble in its domain, then $\lowfnl(1)$ is non-zero     as well, and hence, $\lpr$ is uniquely given by $\lowfnl/\lowfnl(1)$.
  \end{enumerate}
\end{theorem}

\begin{proof}
  \eqref{theo:exactness-reduces-to-coherence:simple}. Assume that the functional $\lowfnl$ is   exact, so $\norm{\lowfnl}<+\infty$.  Let's construct a coherent lower prevision $\lpr$   defined on $\Domain\lowfnl$ such that $\lowfnl=\norm{\lowfnl}\lpr$.  The result is trivial   if $\norm{\lowfnl}=0$, because this holds if and only if $\lowfnl=0$.  Let us assume then   that $\norm{\lowfnl}>0$. The natural extension $\lnex_\lowfnl$ is an exact extension of   $\lowfnl$, and $\norm{\lowfnl}=\norm{\lnex_\lowfnl}=\lnex_\lowfnl(1)$, since $1$ belongs to   the domain $\gambles$ of the exact functional $\lnex_\lowfnl$. Define the functional $\alpr$   on $\gambles$ by $\alpr:=\lnex_\lowfnl/\norm{\lowfnl}=\lnex_\lowfnl/\lnex_\lowfnl(1)$.   Since the exact functional $\lnex_\lowfnl$ is super-additive and positively homogenous, so   is $\alpr$. Moreover, for any gamble $f$ we have that $f\geq\inf f$, so it follows from the   monotonicity of $\lnex_\lowfnl$ and property (E5) that   $\lnex_\lowfnl(f)\geq\lnex_\lowfnl(\inf f)=\lnex_\lowfnl(1)\,\inf f$, whence   $\alpr(f)\geq\inf f$. This tells us that $\alpr$ is a coherent lower prevision on   $\gambles$. Let $\lpr$ be the restriction of $\alpr$ to $\Domain\lowfnl$; since $\lpr$ is   the restriction of a coherent lower prevision, $\lpr$ must be a coherent lower prevision as   well. It follows that for any gamble $f$ in $\Domain\lpr=\Domain\lowfnl$:
  \begin{equation*}
    \lowfnl(f)
    =\lnex_\lowfnl(f)
    =\norm{\lowfnl}\alpr(f)
    =\norm{\lowfnl}\lpr(f),
  \end{equation*}
  whence indeed $\lowfnl=\norm{\lowfnl}\lpr$.

  Let's now proceed to prove that also $\lnex_\lowfnl=\norm{\lowfnl}\lnex_\lpr$. For every   gamble $f$ on $\pspace$, $\lnex_\lowfnl(f)$ is equal to
  \begin{equation*}
    \sup\set{\sum_{k=1}^n\lambda_k\lowfnl(f_k)+\lambda\norm{\lowfnl}}
    {n\in\SetN,\lambda_k\in\SetR_+,\lambda\in\SetR,f_k\in\Domain\lowfnl,
      \sum_{k=1}^n\lambda_kf_k+\lambda\leq f}
  \end{equation*}
  and since $\lowfnl=\norm{\lowfnl}\lpr$, this is equal to
  \begin{equation*}
    \norm{\lowfnl}\sup\set{\sum_{k=1}^n\lambda_k\lpr(f_k)+\lambda}
    {n\in\SetN,\lambda_k\in\SetR_+,\lambda\in\SetR,f_k\in\Domain\lpr,
      \sum_{k=1}^n\lambda_kf_k+\lambda\leq f}
  \end{equation*}
  and therefore, by Eq.~\eqref{eq:nextension}, equal to $\norm{\lowfnl}\lnex_\lpr(f)$. This   completes the proof of the first statement.

  \eqref{theo:exactness-reduces-to-coherence:complex}. If $\lpr$ is a coherent lower   prevision, then it is an exact functional, and therefore so is $\lambda\lpr$ for any   $\lambda\geq0$. Conversely, if $\lowfnl$ is an exact functional, then, by   \eqref{theo:exactness-reduces-to-coherence:simple}, there is a $\lambda$, namely   $\lambda=\norm{\lowfnl}$, and a coherent lower prevision $\lpr$, such that   $\lowfnl=\lambda\lpr$.

  Obviously, whenever the equality $\lowfnl=\lambda\lpr$ holds, for some exact functional   $\lowfnl$, non-negative real $\lambda$, and coherent lower prevision $\lpr$, it also holds   that $\lambda\lnex_\lpr$ is an exact extension of $\lowfnl$ with norm   $\norm{\lambda\lnex_\lpr}=\lambda\lnex_\lpr(1)=\lambda$.

  Moreover, if $1$ belongs to the domain of $\lowfnl$, then   $\lowfnl(1)=\lambda\lpr(1)=\lambda$, so $\lambda$ is uniquely given by $\lowfnl(1)$.   Finally, if also $\lowfnl$ is non-zero for at least one gamble $f$, then, by
  \begin{equation*}
    \lowfnl(1)\,\inf f
    \le\lowfnl(f)
    \le\lowfnl(1)\,\sup f,
  \end{equation*}
  it can only happen that $\lowfnl(1)$ is non-zero as well. Therefore, $\lpr$ is uniquely   given by $\lowfnl/\lowfnl(1)$.
\end{proof}

\begin{corollary}
  A functional whose domain contains at least the constant gamble $1$, is a coherent lower   prevision if and only if it is exact and has norm one.
\end{corollary}

So, the set of exact functionals is the convex cone generated by the set of coherent lower previsions, and natural extension commutes with taking non-negative multiples in the following sense: for any coherent lower prevision $\lpr$ and any non-negative real number $\lambda$, the diagram
\begin{equation*}
  \begin{CD}
    \lpr @>\times\lambda>> \lambda\lpr \\
    @V\text{natural extension}VV @VV\text{natural extension}V \\
    \lnex_\lpr @>\times\lambda>> \lambda\lnex_\lpr
  \end{CD}
\end{equation*}
commutes. In summary, Theorem~\ref{theo:exactness-reduces-to-coherence} establishes a one-to-one and onto correspondence between non-zero exact functionals whose domain contains at least the constant gamble $1$, and pairs $(\lambda,\lpr)$ with $\lambda\in\SetR_+$ and $\lpr$ a coherent lower prevision whose domain contains at least the constant gamble $1$; natural extension is compatible with this correspondence.

When the constant gamble $1$ does not belong to the domain of an exact functional $\lowfnl$, the non-negative real number and coherent lower prevision in Theorem~\ref{theo:exactness-reduces-to-coherence} may not be unique, because $\lowfnl$ may have different exact extensions with different norms. Let's demonstrate this with an example.

\begin{example}
  Let $A$ be any proper subset of\/ $\pspace$, so $A\neq\emptyset$ and $A\neq\pspace$. For any   $\alpha\in(0,1]$, define the coherent lower prevision $\lpr_\alpha$ on the singleton   $\{I_A\}$ by $\lpr_\alpha(I_A):=\alpha$.  Then, clearly, for any $\beta\in(0,1]$,
  \begin{equation*}
    \lpr_\alpha=\frac{\alpha}{\beta}\lpr_\beta,
  \end{equation*}
  and hence, the exact functional $\lpr_\alpha$ can be written in many ways as the product of   a non-negative real number and a coherent lower prevision.

  This also yields an instance of a coherent lower prevision whose norm is different from one,   because when $\alpha\in(0,1)$:
  \begin{align*}
    \norm{\lpr_\alpha} &=\inf\set{c\in\SetR_+} {(\forall\lambda\ge0)(\forall\mu\in\SetR)
      (I_A\geq\lambda I_A+\mu\implies\alpha\geq\lambda\alpha+\mu c)}\\
    &=\inf\{c\in\SetR_+\colon (\forall\lambda\ge 0)(\forall\mu\in\SetR) \\* & \qquad\qquad     \big((0\geq\mu\text{ and }1-\lambda\geq\mu)
    \implies(1-\lambda)\alpha\geq\mu c\big)\}\\
    \intertext{note that the case $\mu=0$ yields $1-\lambda\ge0\implies(1-\lambda)\alpha\ge       0$, which is always satisfied, so} &=\inf\set{c\in\SetR_+}     {(\forall\lambda\ge0)(\forall\mu<0) (\tfrac{1-\lambda}{\mu}\leq 1
      \implies\tfrac{1-\lambda}{\mu}\alpha\leq c\big)}\\
    &=\inf\set{c\in\SetR_+}{(\forall\kappa\le 1)(\kappa\alpha\leq c)}\\
    &=\alpha.
  \end{align*}

  Finally, note that $\lpr_\alpha$ has many exact extensions with different norms: for any   $\beta\ge\alpha$, the functional $\lowfnl_{\alpha,\beta}$ defined by   $\lowfnl_{\alpha,\beta}(I_A):=\alpha$ and $\lowfnl_{\alpha,\beta}(1):=\beta$ is an exact   extension of $\lpr_\alpha$ with norm $\norm{\lowfnl_{\alpha,\beta}}=\beta$.
\end{example}

Theorem~\ref{theo:exactness-reduces-to-coherence} allows us to extend many results for coherent lower previsions to exact functionals, in a straightforward manner. In particular, assume that $\lowfnl$ is an exact functional and that all the relevant gambles below are in the domain of $\lowfnl$, then $\norm{\lowfnl}\inf f\leq\lowfnl(f)\leq\uppfnl(f)\leq\norm{\lowfnl}\sup f$ and consequently $0\leq\lowfnl(\Abs{f})\leq\uppfnl(\Abs{f})\leq\norm{\lowfnl}\sup\Abs{f}$.  Also, both $\Abs{\lowfnl(f)-\lowfnl(g)}\leq\uppfnl(\Abs{f-g})$ and $\Abs{\uppfnl(f)-\uppfnl(g)}\leq\uppfnl(\Abs{f-g})$, and therefore, if a sequence $f_n$ of gambles converges uniformly to a gamble $f$, \textit{i.e.}, if $\sup\Abs{f_n-f}\to0$, then $\lowfnl(f_n)\to\lowfnl(f)$ and $\uppfnl(f_n)\to\uppfnl(f)$: any exact functional, and its conjugate, are (in fact, uniformly) continuous with respect to the supremum norm.

\section{$n$-Monotone functionals}
\label{sec:definition}

We are now ready to start our study of the notion of $n$-monotonicity for (exact) functionals.

A subset $\lattice$ of $\gambles$ is called a \emph{lattice} if it is closed under point-wise maximum $\vee$ and point-wise minimum $\wedge$, \textit{i.e.}, if for all $f$ and $g$ in $\lattice$, both $f\join g$ and $f\meet g$ also belong to $\lattice$.  For instance, the set $\gambles$ of all gambles on $\pspace$ is a lattice.

The following definition is a special case of Choquet's general definition of $n$-monotoni\-city \cite{choquet1953} for functions from an Abelian semi-group to an Abelian group.

\begin{definition}\label{def:n-monotonicity}
  Let $n\in\SetNinfty$, and let\/ $\lowfnl$ be a functional whose domain $\Domain\lowfnl$ is a   lattice of gambles on $\pspace$. Then we call\/ $\lowfnl$ \emph{$n$-monotone} if for all   $p\in\SetN$, $p\leq n$, and all $f$, $f_1$, \dots, $f_p$ in $\Domain\lowfnl$:
  \begin{equation*}
    \sum_{I\subseteq\{1,\dots,p\}}(-1)^\card{I}
    \lowfnl\left(f\meet\bigwedge_{i\in I}f_i\right)\geq0.
  \end{equation*}
  The conjugate of an $n$-monotone functional is called \emph{$n$-alternating}.  An   $\infty$-monotone functional (i.e, a functional which is $n$-monotone for all $n \in \SetN$)   is also called \emph{completely monotone}, and its conjugate \emph{completely alternating}.
\end{definition}
\noindent
In this definition, and further on, we use the convention that for $I=\emptyset$, $\bigwedge_{i\in I}f_i$ simply drops out of the expressions (we could let it be equal to $+\infty$). Clearly, if a functional $\lowfnl$ is $n$-monotone, it is also $p$-monotone for $1\leq p\leq n$. The following proposition gives an immediate alternative characterisation for the $n$-monotonicity for functionals.

\begin{proposition}\label{prop:n-monotonicity-gambles}
  Let $n\in\SetNinfty$, and consider a functional\/ $\lowfnl$ whose domain $\Domain\lowfnl$ is   a lattice of gambles on $\pspace$. Then\/ $\lowfnl$ is $n$-monotone if and only if
  \begin{enumerate}[(i)]
  \item $\lowfnl$ is monotone, \textit{i.e.}, for all $f$ and $g$ in $\Domain\lowfnl$ such     that $f\leq g$, we have\/ $\lowfnl(f)\leq\lowfnl(g)$; and
  \item for all $p\in\SetN$, $2\leq p\leq n$, and all $f_1$, \dots, $f_p$ in $\Domain\lowfnl$:
    \begin{equation*}
      \lowfnl\left(\bigvee_{i=1}^pf_i\right)
      \geq\sum_{\emptyset\neq I\subseteq\{1,\dots,p\}}(-1)^{\card{I}+1}
      \lowfnl\left(\bigwedge_{i\in I}f_i\right).
    \end{equation*}
  \end{enumerate}
\end{proposition}

Exactness guarantees $n$-monotonicity only if $n=1$: any exact functional on a lattice of gambles is monotone 
but not necessarily $2$-monotone, as the following counterexample shows.

\begin{counter}\label{count:coherent-not-2monotone}
  Let $\pspace=\{a,b,c\}$, and consider the lower prevision $\lpr$ defined on $\{1,f\}$ by   $\lpr(f)=\lpr(1)=1$, where $f(a)=0$, $f(b)=1$, $f(c)=2$.  The natural extension $\lnex_\lpr$   of $\lpr$, defined on the set $\gambles$ of all gambles on $\pspace$ (obviously a lattice),   is
  \begin{equation*}
    \lnex_\lpr(g)=\min\left\{g(b),g(c),\frac{g(a)+g(c)}{2}\right\}
  \end{equation*}
  for all gambles $g$ on $\pspace$.  The restriction of $\lnex_\lpr$ to the lattice of   $\{0,1\}$-valued gambles (\textit{i.e.}, indicators) on $\pspace$, is a $2$-monotone   coherent lower probability, simply because any coherent lower probability on a three-element   space is easily seen to be $2$-monotone (see also Walley \cite[p.~58]{walley1981}).   However, $\lnex_\lpr$ is \emph{not} $2$-monotone: $1=\lnex_\lpr(f\join   1)<\lnex_\lpr(f)+\lnex_\lpr(1)-\lnex_\lpr(f\meet 1) =1+1-0.5$, which violates the condition   for 2-monotonicity.
\end{counter}

\begin{theorem}\label{theo:linprevs-completely-everything}
  A linear exact functional\/ $\linfnl$ defined on a lattice of gambles is always completely   monotone and completely alternating.
\end{theorem}
\begin{proof}
  By definition, the linear exact functional $\linfnl$ is the restriction of some linear exact   functional $\alinfnl$ on $\gambles$.  Now recall that $\alinfnl$ is a positive real-valued   linear functional, and apply it to both sides of the following well-known identity (for   indicators of events this is known as the \emph{sieve formula}, or \emph{inclusion-exclusion     principle}, see \cite{aigner1977})
  \begin{equation*}
    \bigvee_{i=1}^pf_i
    =\sum_{\emptyset\neq I\subseteq\{1,\dots,p\}}(-1)^{\card{I}+1}
    \bigwedge_{i\in I}f_i.
  \end{equation*}
  to get
  \begin{equation*}
    \alinfnl\left(\bigvee_{i=1}^pf_i\right)
    =\sum_{\emptyset\neq I\subseteq\{1,\dots,p\}}(-1)^{\card{I}+1}
    \alinfnl\left(\bigwedge_{i\in I}f_i\right).
  \end{equation*}
  Since $\alinfnl$ is also ($1$-)monotone, we derive from   Proposition~\ref{prop:n-monotonicity-gambles} that it is completely monotone, and because in   this case condition (ii) in Proposition~\ref{prop:n-monotonicity-gambles} holds with   equality, it is completely alternating as well. Now recall that $\alinfnl$ and $\linfnl$   coincide on the lattice of gambles $\Domain\linfnl$, that contains all the suprema and   infima in the above expression as soon as the $f_i$ belong to $\Domain\linfnl$.
\end{proof}

The following lemma tells us how to construct $n$-monotone functionals via $\wedge$-homo\-morphisms, and can also be useful for instance to prove that a functional is $n$-monotone, by writing it as a concatenation of a simpler $n$-monotone functional and a $\wedge$-homomorphism. This generalises a similar result by Choquet \cite[Chap.~V, Sect.~23.2, p.~197, and Sect.~24.3, p.~198]{choquet1953} from events (using $\cap$-homomorphisms) to gambles.

A $\wedge$-homomorphism $r$ is a mapping from a lattice to a lattice which preserves the $\wedge$ operation: $r(f\wedge g)=r(f)\wedge r(g)$ for all $f$ and $g$ in the domain of $r$. Note that a $\wedge$-homomorphism is necessarily monotone: $f\ge g$ implies $r(f)\ge r(g)$ [if $f\ge g$, then $f\wedge g=g$, so $r(g)=r(f\wedge g)=r(f)\wedge r(g)$ which can only hold if $r(f)\ge r(g)$].

\begin{lemma}\label{lem:wedgehomomorphisms:nmonotonicity}
  Let $n\in\SetNinfty$, let\/ $\lowfnl$ be an $n$-monotone functional defined on a lattice of   gambles, and let $r$ be a $\wedge$-homomorphism from a lattice of gambles $\Domain r$ to the   lattice of gambles $\Domain\lowfnl$.  Then $\alowfnl:=\lowfnl\circ r$ is an $n$-monotone   functional on $\Domain r$.
\end{lemma}

\begin{proof}
  We prove that the conditions of Proposition~\ref{prop:n-monotonicity-gambles} are satisfied.

  It is easily shown that $\alowfnl$ is monotone, \textit{i.e.}, $\alowfnl(f)\ge\alowfnl(g)$   whenever $f\ge g$ for $f$ and $g$ in $\Domain r$ [use the monotonicity of $r$ and   $\lowfnl$].

  Now, for any $p\in\SetN$, $2\le p\le n$, and any $f_1$, \dots, $f_p\in\Domain r$, it holds   that
  \begin{align*}
    \sum_{\emptyset\neq I\subseteq\{1,\dots,p\}}(-1)^{\Abs{I}+1} \alowfnl\left(\bigwedge_{i\in         I}f_i\right) &=\sum_{\emptyset\neq I\subseteq\{1,\dots,p\}}(-1)^{\Abs{I}+1}
    \lowfnl\left(r\left(\bigwedge_{i\in I}f_i\right)\right)\\
    \intertext{and, since $r$ is a $\wedge$-homomorphism,} &=\sum_{\emptyset\neq       I\subseteq\{1,\dots,p\}}(-1)^{\Abs{I}+1}
    \lowfnl\left(\bigwedge_{i\in I}r(f_i)\right)\\
    \intertext{and since $\lowfnl$ is $n$-monotone,}
    &\le\lowfnl\left(\bigvee_{i=1}^p r(f_i)\right)\\
    \intertext{and, since a $\wedge$-homomorphism is monotone, it holds that $r(f_j)\le       r(\bigvee_{i=1}^p f_i)$ for all $j\in\{1,\dots,p\}$, and hence, $\bigvee_{i=1}^p       r(f_i)\le r(\bigvee_{i=1}^p f_i)$. So, again since $\lowfnl$ is monotone,}     &\leq\lowfnl\left(r\left(\bigvee_{i=1}^pf_i\right)\right) =\alowfnl\left(\bigvee_{i=1}^p       f_i\right).
  \end{align*}
  This establishes the lemma.
\end{proof}

\section{$n$-Monotone set functions}
\label{sec:lower-probabilities}
\subsection{Exactness, natural extension to events, and the inner set function}
\label{sec:lower-probabilities:natxttoevents}
If a lattice of gambles contains only (indicators of) events, we call it a \emph{lattice of   events}. A lattice of events is therefore a collection of subsets of $\pspace$ that is closed under (finite) intersection and union. If it is also closed under set complementation and contains the empty set $\emptyset$, we call it a \emph{field}.

We call \emph{set function} any functional $\lowfnl$ defined on a collection of (indicators of) events.  An $n$-monotone functional on a lattice of events is called an \emph{$n$-monotone   set function}.  A \emph{completely monotone set function} is one that is $\infty$-monotone, or equivalently, $p$-monotone for all $p\in\SetN$.

Let us first study the relationship between $n$-monotonicity and exactness for set functions. Recall that $1$-monotonicity is necessary, but not sufficient, for exactness.  We show in what follows that for $n\geq2$, $n$-monotonicity is sufficient, but not necessary, for exactness. To this end, we consider the \emph{inner set function} $\lowfnl_*$ associated with a monotone set function $\lowfnl$ whose domain $\Domain\lowfnl$ is a lattice of events containing $\emptyset$.  $\lowfnl_*$ is defined by
\begin{equation*}
  \lowfnl_*(A)=\sup\set{\lowfnl(B)}{B\in\Domain\lpr\text{ and }B\subseteq A},
\end{equation*}
for any $A\subseteq\pspace$. Clearly $\lowfnl_*$ is monotone as well, and coincides with $\lowfnl$ on its domain $\Domain\lowfnl$.  But $\lowfnl_*$ is not necessarily real-valued; however, it is real-valued when ($\emptyset$ and) $\pspace$ belong to $\Domain\lowfnl$.

Let's first mention some important known results for $2$-monotone set functions, or lower probabilities (recall that any $n$-monotone set function, for $n\ge 2$, is also $2$-monotone). Note that a coherent lower probability $\lpr$ defined on a lattice of events is $2$-monotone if and only if for all $A$ and $B$ in $\Domain\lpr$:
\begin{equation*}
  \lpr(A\cup B)+\lpr(A\cap B)\geq\lpr(A)+\lpr(B).
\end{equation*}
Walley has shown that a $2$-monotone lower probability $\lpr$ defined \emph{on a field} is coherent if and only if $\lpr(\emptyset)=0$ and $\lpr(\pspace)=1$ (this is a consequence of Walley \cite[Thm.~6.1, p.~55--56]{walley1981}). He has also shown that if $\lpr$ is a coherent $2$-monotone lower probability \emph{on a field}, then its inner set function $\lpr_*$ is $2$-monotone as well and agrees with the natural extension $\lnex_\lpr$ of $\lpr$ on events (see Walley \cite[Thm.~3.1.5, p.~125]{walley1991}).  Applying Theorem~\ref{theo:exactness-reduces-to-coherence}, we get the following result, which summarises Walley's findings and extends them to exact set functions.

\begin{proposition}\label{prop:what-is-known}
  A $2$-monotone set function $\lowfnl$ defined on a field of events is exact if and only if\/   $\lowfnl(\emptyset)=0$. In that case its inner set function $\lowfnl_*$ is $2$-monotone as   well and agrees with the natural extension $\lnex_\lowfnl$ on events.
\end{proposition}

In this section, we generalise these results to $n$-monotone set functions defined on a \emph{lattice} of events containing $\emptyset$ and $\pspace$.

First, we prove that the inner set function preserves $n$-monotonicity; this result is actually due to Choquet \cite[Chapt.~IV, Lem.~18.3]{choquet1953} (once it is noted that Choquet's `interior capacity' coincides with our inner set function).  As the proof in Choquet's paper consists of no more than a hint \cite[p.~186, ll.~6--9]{choquet1953}, we work out the details below.

\begin{theorem}\label{theo:inner-extension-monotone-events}
  Let $n\in\SetNinfty$. Let\/ $\lowfnl$ be a set function defined on a lattice of events   containing $\emptyset$ and $\pspace$. If\/ $\lowfnl$ is $n$-monotone, then its inner set   function $\lowfnl_*$ is $n$-monotone as well.
\end{theorem}
\begin{proof}
  Let $p\in\SetN$, $p\leq n$, and consider arbitrary subsets $B$, $B_1$, \dots, $B_p$ of   $\pspace$.  Fix $\epsilon>0$. Then for each $I\subseteq\{1,\dots,p\}$ it follows from the   definition of $\lowfnl_*$ that there is some $D_I$ in $\Domain\lowfnl$ such that   $D_I\subseteq B\cap\bigcap_{i\in I}B_i$ and
  \begin{equation}\label{eq:inner-extension-intermediate}
    \lowfnl_*\left(B\cap\bigcap_{i\in I}B_i\right)-\epsilon
    \leq\lowfnl(D_I)\leq\lowfnl_*\left(B\cap\bigcap_{i\in I}B_i\right);
  \end{equation}
  note that $\lowfnl_*$ is real-valued since $\emptyset$ and $\pspace$ belong to   $\Domain\lowfnl$. Similarly as before, we use the convention that for $I=\emptyset$, the   corresponding intersection drops out of the expressions (we let it be equal to $\pspace$).   We also let the union of an empty class be equal to $\emptyset$.  Define, for any   $I\subseteq\{1,\dots,p\}$, $E_I=\bigcup_{I\subseteq J\subseteq\{1,\dots,p\}}D_J$, then   clearly $E_I\in\Domain\lowfnl$ and
  \begin{equation*}
    D_I\subseteq E_I\subseteq B\cap\bigcap_{i\in I}B_i.
  \end{equation*}
  Now let $F=E_\emptyset$ and $F_k=E_{\{k\}}\subseteq F$ for $k=1,\dots,p$.  Then $F$ and all   the $F_k$ belong to $\Domain\lowfnl$, and we have for any $K\subseteq\{1,\dots,p\}$ and any   $k\in K$ that $E_K\subseteq E_{\{k\}}=F_k\subseteq B\cap B_k$, whence
  \begin{equation*}
    E_K
    \subseteq\bigcap_{k\in K}F_k
    =F\cap\bigcap_{k\in K}F_k
    \subseteq B\cap\bigcap_{k\in K}B_k.
  \end{equation*}
  Summarising, we find that for every given $\epsilon>0$, there are $F$ and $F_k$ in   $\Domain\lowfnl$, such that for all $I\subseteq\{1,\dots,p\}$
  \begin{equation}\label{eq:inner-extension-intermediate-too}
    D_I\subseteq F\cap\bigcap_{i\in I}F_i\subseteq B\cap\bigcap_{i\in I}B_i
  \end{equation}
  and, using the monotonicity of $\lowfnl_*$ and the fact that it coincides with $\lowfnl$ on   its domain $\Domain\lowfnl$, since $\lowfnl$ is monotone, we deduce from   Eqs.~\eqref{eq:inner-extension-intermediate} and~\eqref{eq:inner-extension-intermediate-too}   that
  \begin{equation*}
    \lowfnl_*\left(B\cap\bigcap_{i\in I}B_i\right)-\epsilon
    \leq\lowfnl\left(F\cap\bigcap_{i\in I}F_i\right)
    \leq\lowfnl_*\left(B\cap\bigcap_{i\in I}B_i\right).
  \end{equation*}
  Consequently, for every $\epsilon>0$ we find that
  \begin{align*}
    \sum_{I\subseteq\{1,\dots,p\}}&(-1)^\card{I} \lowfnl_*\left(B\cap\bigcap_{i\in         I}B_i\right) \\ &=\sum_{\substack{I\subseteq\{1,\dots,p\}\\\text{$I$ even}}}     \lowfnl_*\left(B\cap\bigcap_{i\in I}B_i\right)     -\sum_{\substack{I\subseteq\{1,\dots,p\}\\\text{$I$ odd}}}
    \lowfnl_*\left(B\cap\bigcap_{i\in I}B_i\right)\\
    &\geq\sum_{\substack{I\subseteq\{1,\dots,p\}\\\text{$I$ even}}}     \lowfnl\left(F\cap\bigcap_{i\in I}F_i\right)     -\sum_{\substack{I\subseteq\{1,\dots,p\}\\\text{$I$ odd}}}
    \left[\lowfnl\left(F\cap\bigcap_{i\in I}F_i\right)+\epsilon\right]\\
    &=\sum_{I\subseteq\{1,\dots,p\}}(-1)^\card{I} \lowfnl\left(F\cap\bigcap_{i\in         I}F_i\right)-N_p\epsilon \geq-N_p\epsilon,
  \end{align*}
  where $N_p=2^{p-1}$ is the number of subsets of $\{1,\dots,p\}$ with an odd number of   elements, and the last inequality follows from the $n$-monotonicity of the set function   $\lowfnl$. Since this holds for all $\epsilon>0$, we find that the inner set function   $\lowfnl_*$ is $n$-monotone on the lattice of events $\wp(\pspace)$.
\end{proof}

Recall from Section~\ref{sec:exactfunctionals} that an exact set function on a lattice of events is always monotone, or in other words, $1$-monotone.  In Counterexample~\ref{count:coherent-not-2monotone}, we showed that an exact functional that is $2$-monotone on all events need not be $2$-monotone on all gambles. But at the same time, a set function defined on a field of events can be exact without necessarily being $2$-monotone, as Walley shows (for the special case of coherent lower probabilities) in \cite[p.~51]{walley1981}.  Conversely, a $2$-monotone set function defined on a lattice of events need not be exact: it suffices to consider any constant non-zero set function on $\wp(\pspace)$. Below, we give simple necessary and sufficient conditions for the exactness of an $n$-monotone set function, we characterise its natural extension, and we prove that the natural extension of an $n$-monotone exact set function to all events is again an $n$-monotone exact set function.

\begin{proposition}\label{prop:charac-coherence}
  Let\/ $\lowfnl$ be an $n$-monotone set function ($n\in\SetNinfty$, $n\ge 2$) defined on a   lattice of events that contains $\emptyset$ and $\pspace$.  Then $\lowfnl$ is exact if and   only if\/ $\lowfnl(\emptyset)=0$.
\end{proposition}

\begin{proof}
  Clearly, $\lowfnl(\emptyset)=0$ is necessary for exactness.  Conversely, by   Theorem~\ref{theo:inner-extension-monotone-events}, the inner set function $\lowfnl_*$ of   $\lowfnl$ to all events is also $n$-monotone, and hence $2$-monotone. Now, $\lowfnl_*$ is   defined on a field, so $\lowfnl_*$ must be exact as we already argued before (see   Proposition~\ref{prop:what-is-known}, or alternatively, apply   Theorem~\ref{theo:exactness-reduces-to-coherence} and Walley \cite[Thm.~6.1,   p.~55--56]{walley1981}).  Consequently $\lowfnl$ is exact as well.
\end{proof}

The following proposition relates the natural extension $\lnex_\lowfnl$ of an exact $n$-monotone set function $\lowfnl$ with the inner set function $\lowfnl_*$.

\begin{proposition}\label{pr:inner-equal-natural}
  Let\/ $\lowfnl$ be an exact $n$-monotone set function ($n\in\SetNinfty$, $n\ge 2$) defined   on a lattice of events that contains $\emptyset$ and $\pspace$.  Then its natural extension   $\lnex_\lowfnl$ restricted to events is an $n$-monotone exact set function as well, and it   coincides with the inner set function $\lowfnl_*$ of\/ $\lowfnl$.
\end{proposition}
\begin{proof}
  Take $A\subseteq\pspace$. Then, for any $\alinfnl$ in $\solp(\lowfnl)$, since $\alinfnl$ is   monotone and dominates $\lowfnl$,
  \begin{equation*}
    \alinfnl(A)\geq\sup_{\substack{B\subseteq A,B\in\Domain\lowfnl}}\alinfnl(B)
    \geq\sup_{\substack{B\subseteq A,B\in\Domain\lowfnl}}\lowfnl(B)
    =\lowfnl_*(A).
  \end{equation*}
  Since we know that $\lnex_\lowfnl(A)=\min\set{\alinfnl(A)}{\alinfnl\in\solp(\lowfnl)}$, we   deduce that $\lnex_\lowfnl(A)\geq\lowfnl_*(A)$ for all $A\subseteq\pspace$.

  Conversely, from Theorem \ref{theo:inner-extension-monotone-events}, $\lowfnl_*$ is   $n$-monotone if $\lowfnl$ is, and applying Proposition~\ref{prop:charac-coherence},   $\lowfnl_*$ is an exact extension of $\lowfnl$ to all events. Moreover,   $\norm{\lowfnl_*}=\lowfnl_*(\pspace)=\lowfnl(\pspace)=\norm{\lowfnl}$ (see   Theorem~\ref{theo:maass-results} or Maa\ss\ \cite[Prop.~1.2.4]{maass2003}), and therefore   $\lowfnl_*$ must dominate the natural extension $\lnex_\lowfnl$ of $\lowfnl$ (Maa\ss\   \cite[Prop.~1.2.7(a)]{maass2003}), whence also $\lnex_\lowfnl(A)\leq\lowfnl_*(A)$ for all   $A\subseteq\pspace$.
\end{proof}

Note that this shows in particular that the natural extension of an $n$-monotone exact set function to all events is also $n$-monotone. This result will be generalised in the following section.

\subsection{Natural extension to all gambles, and the Choquet integral}
\label{sec:lower-probabilities:natxttogambles}

Walley \cite[p.~56]{walley1981} has shown that the natural extension $\lnex_\lpr$ to all gambles of a coherent $2$-monotone lower probability $\lpr$ defined on the set $\wp(\pspace)$ of all events, is given by the Choquet functional with respect to $\lpr$.  Hence, by Theorem~\ref{theo:exactness-reduces-to-coherence}, the natural extension of an exact $2$-monotone set function $\lowfnl$ on $\wp(\pspace)$ is given by
\begin{equation}\label{eq:nex-2monotone}
  \lnex_\lowfnl(f)
  =(C)\int f\dif\lowfnl
  =\norm{\lowfnl}\inf f+(R)\int_{\inf f}^{\sup f}\ddf{f}{\lowfnl}(x)\dif x,
\end{equation}
where the integral on the right-hand side is a Riemann integral, and the function $\ddf{f}{\lowfnl}$ defined by $\ddf{f}{\lowfnl}(x)=\lowfnl(\{f\geq x\})$, is the \emph{decreasing distribution function of $f$} with respect to $\lowfnl$; note that $\ddf{f}{\lowfnl}$ is always bounded and non-increasing, and therefore always Riemann integrable.  We have used the common notation $\{f\geq x\}$ for the set $\set{\omega\in\pspace}{f(\omega)\geq x}$.

Eq.~\eqref{eq:nex-2monotone} tells us also that $\lnex_\lowfnl$ is \emph{comonotone additive} on $\gambles$, because that is a property of any Choquet functional associated with a monotone set function on a field (see Denneberg \cite[Prop.~5.1]{denneberg1994}): if two gambles $f$ and $g$ are \emph{comonotone} in the sense that
\begin{equation*}
  (\forall\omega_1,\omega_2\in\pspace)
  (f(\omega_1)<f(\omega_2)\implies g(\omega_1)\leq g(\omega_2)),
\end{equation*}
then $\lnex_\lowfnl(f+g)=\lnex_\lowfnl(f)+\lnex_\lowfnl(g)$.

By Proposition~\ref{pr:inner-equal-natural}, we may assume that a $2$-monotone exact set function defined on a lattice of events that contains $\emptyset$ and $\pspace$, is actually defined on all of $\wp(\pspace)$, since we can extend it to $\wp(\pspace)$ using the inner set function (or, natural extension) $\lowfnl_*$, which is still $2$-monotone.  Moreover, the natural extension of $\lowfnl$ to all gambles coincides with the natural extension of $\lowfnl_*$ to all gambles, because of the transitivity property mentioned at the end of Section~\ref{sec:exactfunctionals}. This means that Eq.~\eqref{eq:nex-2monotone} also holds for $2$-monotone exact set functions defined on a lattice of events.  Since any $n$-monotone set function, for $n\ge 2$, is also $2$-monotone, we conclude:

\begin{theorem}\label{theo:natext-choquet-lattice}
  Let $n\in\SetNinfty$, $n\ge 2$, and let\/ $\lowfnl$ be an $n$-monotone exact set function   defined on a lattice of events that contains both $\emptyset$ and $\pspace$. Then its   natural extension $\lnex_\lowfnl$ to the set of all gambles is given by
  \begin{equation*}
    \lnex_\lowfnl(f)=(C)\int f\dif\lowfnl_*
    =\norm{\lowfnl}\inf f
    +(R)\int_{\inf f}^{\sup f}\lowfnl_*(\{f\geq x\})\dif x.
  \end{equation*}
\end{theorem}

We already know from Theorem~\ref{theo:inner-extension-monotone-events} that the natural extension of an $n$-monotone exact set function to the set of all events, is $n$-monotone as well. This result holds also for the natural extension to gambles.

\begin{theorem}\label{theo:natext-monotone-gambles}
  Let $n\in\SetNinfty$, $n\ge 2$, and let\/ $\lowfnl$ be an exact set function, defined on a   lattice of events that contains $\emptyset$ and $\pspace$. If\/ $\lowfnl$ is $n$-monotone,   then its natural extension $\lnex_\lowfnl$ is $n$-monotone as well.
\end{theorem}
\begin{proof}
  Let $p\in\SetN$, $p\leq n$, and let $f$, $f_1$, \dots, $f_p$ be arbitrary gambles on   $\pspace$.  Let
  \begin{equation*}
    a=\min\{\inf f,\min_{k=1}^p\inf f_k\} \text{ , }
    b=\max\{\sup f,\max_{k=1}^p\sup f_k\}.
  \end{equation*}
  Consider $I\subseteq\{1,\dots,p\}$ then $a\leq\inf(f\meet\bigwedge_{i\in I}f_i)$ and   $b\geq\sup(f\meet\bigwedge_{i\in I}f_i)$. It is easily verified that
  \begin{equation*}
    \lnex_\lowfnl\left(f\meet\bigwedge_{i\in I}f_i\right)
    =\norm{\lowfnl}a+(R)\int_a^b\ddf{f\meet\bigwedge_{i\in I}f_i}{\lowfnl_*}(x)\dif x.
  \end{equation*}
  Since it is obvious that for any $x$ in $\SetR$
  \begin{equation*}
    \ddf{f\meet\bigwedge_{i\in I}f_i}{\lowfnl_*}(x)
    =\lowfnl_*\left(\{f\geq x\}\cap\bigcap_{i\in I}\{f_i\geq x\}\right),
  \end{equation*}
  it follows from the $n$-monotonicity of $\lowfnl_*$ (see   Theorem~\ref{theo:inner-extension-monotone-events}) that for all real $x$
  \begin{equation*}
    \sum_{I\subseteq\{1,\dots,p\}}(-1)^\card{I}
    \ddf{f\meet\bigwedge_{i\in I}f_i}{\lowfnl_*}(x)\geq0.
  \end{equation*}
  If we take the Riemann integral over $[a,b]$ on both sides of this inequality, and recall   moreover that $\sum_{I\subseteq\{1,\dots,p\}}(-1)^\card{I}=0$, we get
  \begin{equation*}
    \sum_{I\subseteq\{1,\dots,p\}}(-1)^\card{I}
    \lnex_\lowfnl\left(f\meet\bigwedge_{i\in I}f_i\right)\geq0.
  \end{equation*}
  This tells us that $\lnex_\lowfnl$ is $n$-monotone.
\end{proof}

We deduce in particular from this result that given an $n$-monotone exact set function defined on $\wp(\pspace)$, the functional that we can define on $\gambles$ by means of its Choquet functional is $n$-monotone and exact.  Since trivially the converse also holds, we deduce that the Choquet functional with respect to an exact set function $\lowfnl$ on $\wp(\pspace)$ is $n$-monotone if and only if $\lowfnl$ is. This generalises a result by Walley \cite[Thm.~6.4]{walley1981}.

\begin{corollary}\label{cor:n-mono:equivalence:lpr:lnex:choqint}
  Let\/ $\lowfnl$ be any exact set function defined on a lattice of events containing both   $\emptyset$ and $\pspace$. Let $n\in\SetNinfty$, $n\ge2$.  Then $\lowfnl$ is $n$-monotone,   if and only if\/ $\lnex_\lowfnl$ is $n$-monotone, if and only if\/   $(C)\int\cdot\dif\lowfnl_*$ is $n$-monotone.
\end{corollary}
\begin{proof}
  If $\lowfnl$ is $n$-monotone, then $\lnex_\lowfnl$ is $n$-monotone by   Theorem~\ref{theo:natext-monotone-gambles}.

  If $\lnex_\lowfnl$ is $n$-monotone, then $\lowfnl$ is $n$-monotone since $\lnex_\lowfnl$ is   an extension of $\lowfnl$ (because $\lowfnl$ is exact), and so, by   Theorem~\ref{theo:natext-choquet-lattice}, $\lnex_\lowfnl$ must coincide with   $(C)\int\cdot\dif\lowfnl_*$, which must therefore be $n$-monotone as well.

  Finally, if $(C)\int\cdot\dif\lowfnl_*$ is $n$-monotone, then $\lowfnl_*$ must be   $n$-monotone since $(C)\int\cdot\dif\lowfnl_*$ is an extension of $\lowfnl_*$. But,   $\lowfnl_*$ is also an extension of $\lowfnl$ (because $\lowfnl$ is also $1$-monotone), so,   $\lowfnl$ is $n$-monotone as well. This completes the chain.
\end{proof}

\subsection{Application: minimum preserving functionals are completely monotone}

A functional $\lowfnl$ defined on a lattice of gambles is called \emph{minimum preserving} if $\lowfnl(f\wedge g)=\lowfnl(f)\wedge\lowfnl(g)$ for all $f$ and $g$ in $\Domain\lowfnl$, that is, if it is a $\wedge$-homomorphism between its domain and $\SetR$.

Now, $\wedge$-homomorphisms (Lemma~\ref{lem:wedgehomomorphisms:nmonotonicity}) and natural extension (Theorem~\ref{theo:natext-monotone-gambles}) provide two ways to deduce $n$-monotone functionals from other $n$-monotone functionals.  Combining these results we easily obtain that any minimum preserving functional is completely monotone. This generalises a result by Nguyen \cite[Thm.~1, p.~363--364]{nguyen1997} from set functions to functionals. Also note that, in contradistinction to Nguyen's proof, our proof does not rely on combinatorics.

\begin{theorem}\label{thm:minpreserving:completelymonotone}
  Any minimum preserving functional defined on a lattice of gambles is completely monotone.
\end{theorem}

\begin{proof}
  Let $\lowfnl$ be a minimum preserving functional defined on a lattice of gambles.  Define   the lower probability $\alpr$ on $\{\emptyset,\pspace\}$ by $\alpr(\emptyset)=0$ and   $\alpr(\pspace)=1$.  Clearly, $\alpr$ is a completely monotone exact set function (it is   even a probability charge).  Hence, its natural extension $\nex_\alpr$ to $\gambles$ is   completely monotone, by Theorem~\ref{theo:natext-monotone-gambles}. Since $\alpr$ is   dominated by all linear previsions on $\gambles$ (and in particular by the degenerate   probability distributions on some $\omega \in \pspace$), it's not difficult to see that   $\lnex_\alpr(f)=\inf f$ for all gambles $f$ on $\pspace$.

  Now, define the mapping $r\colon\Domain\lowfnl\to\gambles$ by $r(f)(\omega):=\lowfnl(f)$ for   all $f$ in $\Domain\lowfnl$ and all $\omega\in\pspace$. Since $\lowfnl$ is minimum   preserving, $r$ is a $\wedge$-homomorphism. Observe that $\lowfnl=\lnex_\alpr\circ r$, and   apply Lemma~\ref{lem:wedgehomomorphisms:nmonotonicity}.
\end{proof}

As an example, the \emph{vacuous lower prevision} relative to a non-empty subset $A$ of $\pspace$, given by
\begin{equation*}
  \lpr_A(f):=\inf_{\omega\in A}f(\omega),
\end{equation*}
for all $f$ in $\gambles$, is minimum preserving. So, $\lpr_A$ is an instance of a completely monotone lower prevision on $\gambles$.

\section{Representation results}
\label{sec:representation}

Let us now focus on the notion of $n$-monotonicity we have given for functionals.  If $\lowfnl$ is a monotone functional on a lattice of gambles that contains all constant gambles, then its \emph{inner extension} $\lowfnl_*$ is given by
\begin{equation}\label{eq:inner-extension}
  \lowfnl_*(f)=\sup\set{\lowfnl(g)}{g\in\Domain\lpr\text{ and }g\leq f}.
\end{equation}
for all gambles $f$ on $\pspace$. Clearly this inner extension is monotone as well, and it coincides with $\lowfnl$ on its domain $\Domain\lowfnl$. The following result in some sense generalises Theorem~\ref{theo:inner-extension-monotone-events}.

\begin{theorem}\label{theo:inner-extension-monotone-gambles}
  Let $n\in\SetNinfty$.  Let\/ $\lowfnl$ be a functional defined on a lattice of gambles that   contains all constant gambles. If\/ $\lowfnl$ is $n$-monotone, then $\lowfnl_*$ is   $n$-monotone as well.
\end{theorem}
\begin{proof}
  Let $p\in\SetN$, $p\le n$, and consider arbitrary gambles $f$, $f_1$, \dots, $f_p$ on   $\pspace$.  Fix $\epsilon>0$. Since $\Domain\lowfnl$ is assumed to contain all constant   gambles, and since gambles are bounded, we see that for each $I\subseteq\{1,\dots,p\}$ there   is some $g_I$ in $\Domain\lowfnl$ such that $g_I\leq f\meet\bigwedge_{i\in I}f_i$ and
  \begin{equation*}
    \lowfnl_*\left(f\meet\bigwedge_{i\in I}f_i\right)-\epsilon
    \leq\lowfnl(g_I)\leq\lowfnl_*\left(f\meet\bigwedge_{i\in I}f_i\right).
  \end{equation*}
  Define, for any $I\subseteq\{1,\dots,p\}$, $h_I=\bigvee_{I\subseteq     J\subseteq\{1,\dots,p\}}g_J$, then clearly $h_I\in\Domain\lowfnl$ and
  \begin{equation*}
    g_I\leq h_I\leq f\meet\bigwedge_{i\in I}f_i.
  \end{equation*}
  Now consider the gambles $q=h_\emptyset$ and $q_k=h_{\{k\}}\leq q$ for $k=1,\dots,p$. Then   $q$ and all the $q_k$ belong to $\Domain\lowfnl$, and we have for any   $K\subseteq\{1,\dots,p\}$ and any $k\in K$ that $h_K\leq h_{\{k\}}=q_k\leq f\meet f_k$,   whence
  \begin{equation*}
    h_K\leq\bigwedge_{k\in K}q_k=q\meet\bigwedge_{k\in K}q_k
    \leq f\meet\bigwedge_{k\in K}f_k.
  \end{equation*}
  Summarising, we find that for every given $\epsilon>0$, there are $q$ and $q_k$ in   $\Domain\lowfnl$, such that for all $I\subseteq\{1,\dots,p\}$
  \begin{equation*}
    g_I\leq q\meet\bigwedge_{i\in I}q_i\leq f\meet\bigwedge_{i\in I}f_i
  \end{equation*}
  and, using the monotonicity of $\lowfnl_*$ and the fact that it coincides with $\lowfnl$ on   its domain $\Domain\lowfnl$, since $\lowfnl$ is monotone,
  \begin{equation*}
    \lowfnl_*\left(f\meet\bigwedge_{i\in I}f_i\right)-\epsilon
    \leq\lowfnl\left(q\meet\bigwedge_{i\in I}q_i\right)
    \leq\lowfnl_*\left(f\meet\bigwedge_{i\in I}f_i\right).
  \end{equation*}
  Consequently, for every $\epsilon>0$ we find that
  \begin{align*}
    \sum_{I\subseteq\{1,\dots,p\}}&(-1)^\card{I} \lowfnl_*\left(f\meet\bigwedge_{i\in         I}f_i\right) \\ &=\sum_{\substack{I\subseteq\{1,\dots,p\}\\\text{$I$ even}}}     \lowfnl_*\left(f\meet\bigwedge_{i\in I}f_i\right)     -\sum_{\substack{I\subseteq\{1,\dots,p\}\\\text{$I$ odd}}}
    \lowfnl_*\left(f\meet\bigwedge_{i\in I}f_i\right)\\
    &\geq\sum_{\substack{I\subseteq\{1,\dots,p\}\\\text{$I$ even}}}     \lowfnl\left(q\meet\bigwedge_{i\in I}q_i\right)     -\sum_{\substack{I\subseteq\{1,\dots,p\}\\\text{$I$ odd}}}
    \left[\lowfnl\left(q\meet\bigwedge_{i\in I}q_i\right)+\epsilon\right]\\
    &=\sum_{I\subseteq\{1,\dots,p\}}(-1)^\card{I} \lowfnl\left(q\meet\bigwedge_{i\in         I}q_i\right) -N_p\epsilon \geq-N_p\epsilon,
  \end{align*}
  where $N_p=2^{p-1}$ is the number of subsets of $\{1,\dots,p\}$ with an odd number of   elements, and the last inequality follows from the $n$-monotonicity of $\lowfnl$. Since this   holds for all $\epsilon>0$, we find that $\lowfnl_*$ is $n$-monotone on the lattice of   gambles $\gambles$.
\end{proof}

We now investigate whether a result akin to Theorem~\ref{theo:natext-monotone-gambles} holds for $n$-monotone exact functionals: when will the natural extension of an $n$-monotone exact functional be $n$-monotone? For Theorem~\ref{theo:natext-monotone-gambles}, we needed the domain of the set function to be a lattice of events containing $\emptyset$ and $\pspace$. It turns out that for our generalisation we also have to impose a similar condition on the domain: it will have to be a linear lattice containing all constant gambles. Recall that a subset $\linlat$ of $\gambles$ is called a \emph{linear lattice} if $\linlat$ is a linear space under point-wise addition and scalar multiplication with real numbers, and if it is moreover closed under point-wise minimum $\meet$ and point-wise maximum $\join$.

Consider an exact functional whose domain is a linear lattice of gambles that contains all constant gambles.  Then its natural extension to the set of all gambles $\gambles$ is precisely its inner extension $\lowfnl_*$, by Walley \cite[Thm.~3.1.4]{walley1991} and Theorem~\ref{theo:exactness-reduces-to-coherence}. This leads at once to the following theorem, which is a counterpart of Theorem~\ref{theo:natext-monotone-gambles} for $n$-monotone exact functionals.

\begin{theorem}\label{theo:natext-monotone-gambles-linlat}
  Let $n\in\SetNinfty$, and let\/ $\lowfnl$ be an exact functional defined on a linear lattice   of gambles that contains all constant gambles. If\/ $\lowfnl$ is $n$-monotone, then its   natural extension $\lnex_\lowfnl$ is equal to its inner extension $\lowfnl_*$, and is   therefore $n$-monotone as well.
\end{theorem}
\noindent
Counterexample~\ref{count:coherent-not-2monotone} tells us that this result cannot be extended to lattices of gambles that are not at the same time linear spaces.

We have not made any mention yet of the Choquet integral in relation to the natural extension. It turns out that, to some extent, there is also a relationship between both concepts. Consider a linear lattice of gambles $\linlat$ that contains all constant gambles. Then the set
\begin{equation*}
  \field_\linlat=\set{A\subseteq\pspace}{I_A\in\linlat}
\end{equation*}
of events that belong to $\linlat$ is a field of subsets of $\pspace$.  Let us denote by $\gambles_{\field_\linlat}$ the uniformly closed linear lattice
\begin{equation*}
  \gambles_{\field_\linlat}=\cl(\linspan(I_{\field_\linlat})),
\end{equation*}
where $I_{\field_\linlat}=\set{I_A}{I_A\in\linlat}$, $\cl(\cdot)$ denotes uniform closure, and $\linspan(\cdot)$ takes the linear span.  Observe that $\gambles_{\field_\linlat}$ contains all constant gambles as well. We call its elements \emph{$\field_\linlat$-measurable   gambles}. Every $\field_\linlat$-measurable gamble is a uniform limit of \emph{$\field_\linlat$-simple gambles}, \textit{i.e.}, elements of $\linspan(I_{\field_\linlat})$. Moreover, $\gambles_{\field_\linlat}\subseteq\cl(\linlat)$.

\begin{theorem}\label{theo:monotone-representation}
  Let\/ $\lowfnl$ be an $n$-monotone exact functional on a linear lattice of gambles $\linlat$   that contains all constant gambles.  Then $\lowfnl$ has a unique exact extension to   $\cl(\linlat)$, and this extension is $n$-monotone as well.  Denote by $\alowfnl$ the   restriction of\/ $\lowfnl$ to $\field_\linlat$. Then for all $f$ in   $\gambles_{\field_\linlat}$,
  \begin{equation*}
    \lnex_\lowfnl(f)=\lnex_\alowfnl(f)=(C)\int f\dif\alowfnl_*
    =\norm{\lowfnl}\inf f
    +(R)\int_{\inf f}^{\sup f}\alowfnl_*(\{f\geq x\})\dif x.
  \end{equation*}
  Consequently, $\lnex_\lowfnl$ is both $n$-monotone and comonotone additive on   $\gambles_{\field_\linlat}$.
\end{theorem}

\begin{proof}
  Let us first show that $\lowfnl$ has a unique exact extension to $\cl(\linlat)$. Let   $\aalowfnl$ be any such exact extension.  If we can show that $\aalowfnl$ coincides with   $\lnex_\lowfnl$ on $\cl(\linlat)$, then we have established uniqueness.  Consider any   element $h$ in $\cl(\linlat)$.  Then there is a sequence $g_n$ of gambles in $\linlat$ that   converges uniformly to $h$.  Since both $\aalowfnl$ and $\lnex_\lowfnl$ coincide with   $\lowfnl$ on $\linlat$, and are uniformly continuous on their domain $\cl(\linlat)$, because   they are exact, we indeed find that
  \begin{equation*}
    \aalowfnl(h)=\lim_{n\to\infty}\aalowfnl(g_n)
    =\lim_{n\to\infty}\lnex_\lowfnl(g_n)=\lnex_\lowfnl(h).
  \end{equation*}

  Let's now prove the equalities.  Since $\lowfnl$ is $n$-monotone and exact, its restriction   $\alowfnl$ to the field $\field_\linlat$ is an $n$-monotone exact set function. By   Theorem~\ref{theo:natext-choquet-lattice}, the natural extension $\lnex_\alowfnl$ of   $\alowfnl$ to the set $\gambles$ of all gambles is the Choquet functional associated with   the $n$-monotone inner set function $\alowfnl_*$ of $\alowfnl$: for any gamble $f$ on   $\pspace$,
  \begin{equation*}
    \lnex_\alowfnl(f)
    =(C)\int f\dif\alowfnl_*
    =\norm{\alowfnl}\inf f+(R)\int_{\inf f}^{\sup f}\alowfnl_*(\{f\geq x\})\dif x,
  \end{equation*}
  and note that $\norm{\alowfnl}=\alowfnl(1)=\lowfnl(1)=\norm{\lowfnl}$.

  Finally, to prove that $\lnex_\alowfnl$ and $\lnex_\lowfnl$ coincide on the subset   $\gambles_{\field_\linlat}$ of $\cl(\linlat)$, observe that suffices to prove that   $\lnex_\alowfnl$ and $\lowfnl$ coincide on $\linspan(I_{\field_\linlat})$, since   $\lnex_\alowfnl$ and $\lnex_\lowfnl$ are guaranteed by exactness to be continuous, and since   $\lnex_\lowfnl$ and $\lowfnl$ coincide on $\linspan(I_{\field_\linlat})\subseteq\linlat$,   because $\lowfnl$ is exact on $\linlat$. Let therefore $h$ be any element of   $\linspan(I_{\field_\linlat})$, \textit{i.e.}, let $h$ be an $\field_\linlat$-simple   gamble. Then we can always find a natural number $n\geq1$, real $\mu_1$, real non-negative   $\mu_2$, \dots, $\mu_n$, and nested sets $F_2\supseteq\dots\supseteq F_n$ such that
  \begin{equation*}
    h=\mu_1+\sum_{k=2}^n\mu_kI_{F_k}.
  \end{equation*}
  It then follows from the comonotone additivity of the Choquet integral that
  \begin{equation*}
    \lnex_\alowfnl(h)=\lowfnl(\mu_1)+\sum_{k=2}^n\mu_k\alowfnl(F_k).
  \end{equation*}
  On the other hand, it follows from the exactness and the $2$-monotonicity of $\lowfnl$ that
  \begin{align*}
    \lowfnl(h)
    &=\lowfnl(\mu_1)+\lowfnl\left(\sum_{k=2}^n\mu_kI_{F_k}\right)\\
    &=\lowfnl(\mu_1)-\lowfnl(\mu_2)+\lowfnl\left(\sum_{k=2}^n\mu_kI_{F_k}\right)+\lowfnl(\mu_2)\\
    &\leq\lowfnl(\mu_1)-\lowfnl(\mu_2)+\lowfnl\left(\mu_2\join\sum_{k=2}^n\mu_kI_{F_k}\right)     
    +\lowfnl\left(\mu_2\meet\sum_{k=2}^n\mu_kI_{F_k}\right).
  \end{align*}
  Now it is easily verified that
  \begin{equation*}
    \mu_2\join\sum_{k=2}^n\mu_kI_{F_k}=\mu_2+\sum_{k=3}^n\mu_kI_{F_k}
    \text{ and }
    \mu_2\meet\sum_{k=2}^n\mu_kI_{F_k}=\mu_2I_{F_2},
  \end{equation*}
  and consequently, again using the exactness and the $2$-monotonicity of $\lowfnl$, the fact   that $\alowfnl$ coincides with $\lowfnl$ on $\field_\linlat$, and continuing in the same   fashion,
  \begin{align*}
    \lowfnl(h)     &\leq\lowfnl(\mu_1)-\lowfnl(\mu_2)+\lowfnl\left(\mu_2+\sum_{k=3}^n\mu_kI_{F_k}\right)
    +\lowfnl\left(\mu_2I_{F_2}\right)\\
    &=\lowfnl(\mu_1)+\mu_2\alowfnl(F_2)+\lowfnl\left(\sum_{k=3}^n\mu_kI_{F_k}\right)\\
    &\leq\lowfnl(\mu_1)+\mu_2\alowfnl(F_2)+\mu_3\alowfnl(F_3)
    +\lowfnl\left(\sum_{k=4}^n\mu_kI_{F_k}\right)\\
    &\qquad\vdots\\
    &\leq\lowfnl(\mu_1)+\sum_{k=2}^n\mu_k\alowfnl(F_k).
  \end{align*}
  This tells us that $\lnex_\alowfnl(h)\geq\lowfnl(h)$. On the other hand, since $\lowfnl$ is   an exact extension of $\alowfnl$ with the same norm, and since the natural extension   $\lnex_\alowfnl$ is the point-wise smallest exact extension of $\alowfnl$ with the same   norm, we also find that $\lnex_\alowfnl(h)\leq\lowfnl(h)$. This tells us that $\lowfnl$ and   $\lnex_\alowfnl$ indeed coincide on $\linspan(I_{\field_\linlat})$.
\end{proof}

Walley has shown in \cite{walley1991} that in general coherent lower previsions (and hence, exact functionals) are not determined by their values on events. But the preceding theorem tells us that for exact functionals that are $2$-monotone and defined on a sufficiently rich domain, we can somewhat improve upon this negative result: on $\field_\linlat$-measurable gambles, the natural extension $\lnex_\lowfnl$ of an $n$-monotone exact functional $\lowfnl$ is completely determined by the values that $\lowfnl$ assumes on the events in $\field_\linlat$. Nevertheless, the following counterexample tells us that in general, we cannot expect to take this result beyond the set $\gambles_{\field_\linlat}$ of $\field_\linlat$-measurable gambles.

\begin{counter}\label{count:not-beyond-measurable}
  Let $\pspace$ be the closed unit interval $[0,1]$ in $\SetR$, and let $\lpr$ be the lower   prevision on the lattice $\linlat$ of all continuous gambles on $\pspace$, defined by   $\lpr(f)=f(0)$ for any $f$ in $\linlat$.  Since $\lpr$ is actually a linear prevision, it   must be completely monotone (see Theorem~\ref{theo:linprevs-completely-everything}).   Observe that $\linlat$ is a uniformly closed linear lattice that contains all constant   gambles.  Moreover, $\field_\linlat=\{\emptyset,\pspace\}$, so $\gambles_{\field_\linlat}$   is the set of all constant gambles, and the natural extension $\lnex_\alpr$ of the   restriction $\alpr$ of $\lpr$ to $\field_\linlat$ is the vacuous lower prevision on   $\gambles$: $\lnex_\alpr(f)=\inf f$ for all gambles $f$ on $\pspace$. Therefore, for any $g$   in $\linlat$ such that $g(0)>\inf g$, it follows that $\lnex_\alpr(g)<\lpr(g)$: the equality   in Theorem~\ref{theo:monotone-representation} holds only for those gambles in $\linlat$ that   satisfy $g(0)=\inf g$.
\end{counter}

So we conclude that, in general, an $n$-monotone exact functional $\lowfnl$ defined on a linear lattice of gambles that contains the constant gambles, cannot be written (on its entire domain) as a Choquet functional associated with its restriction $\alowfnl$ to events.

Instead, however, we can represent such $n$-monotone exact functionals by a Choquet integral with respect to the restriction to events of their inner extension, and this Choquet integral also immediately provides us with an alternative expression for the natural extension. This is because $2$-monotonicity and comonotone additivity are equivalent under exactness.

\begin{theorem}\label{theo:equivalence}
  Let\/ $\lowfnl$ be an exact functional defined on a linear lattice of gambles that contains   all constant gambles.  Then $\lowfnl$ is comonotone additive if and only if it is   $2$-monotone, and in both cases we have for all $f$ in $\Domain\lowfnl$
  \begin{equation*}
    \lowfnl(f)=(C)\int f\dif\lowfnl_*
    =\norm{\lowfnl}\inf f
    +(R)\int_{\inf f}^{\sup f}\lowfnl_*(\{f\geq x\})\dif x.
  \end{equation*}
\end{theorem}
\begin{proof}
  Let us first prove the direct implication. Assume that $\lowfnl$ is comonotone additive. Let   us define $\linlat_+:=\set{f\in\Domain\lowfnl}{f\geq0}$, and let $\lowfnl_+$ be the   restriction of $\lowfnl$ to $\linlat_+$. This functional is also exact and comonotone   additive, and it is defined on a class of non-negative gambles.  Moreover, given $f$ in   $\linlat_+$ and $a\geq0$, the gambles $af$, $f\wedge a$ and $f-f\wedge a$ belong to   $\linlat_+$ because $\Domain\lowfnl$ is a linear lattice that contains the constant gambles   and all the above gambles are trivially non-negative. Hence, we may apply Greco's   representation theorem (see \cite[Thm.~13.2]{denneberg1994}; the conditions (iv) and (v)   there are trivially satisfied because all elements in $\linlat_+$ are bounded), and conclude   that there is a monotone set function $\mu$ on $\wp(\pspace)$ with $\mu(\emptyset)=0$ and   $\mu(\pspace)=\lowfnl_+(1)=\lowfnl(1)=\norm{\lowfnl}$ such that for all $f$ in $\linlat_+$:
  \begin{equation*}
    \lowfnl_+(f)=(C)\int f\dif\mu.
  \end{equation*}
  Consider now any $f$ in $\Domain\lowfnl$. Since $f$ is bounded, and exactness implies that   $\lowfnl(f+a)=\lowfnl(f)+\norm{\lowfnl}a$ for all $a$ in $\SetR$, this also implies that   $\norm{\lowfnl}\inf f+\lowfnl_+(f-\inf f)=\lowfnl(f)$, whence
  \begin{equation}\label{eq:equality-everywhere}
    \lowfnl(f)=\norm{\lowfnl}\inf f+(C)\int[f-\inf f]\dif\mu
    =(C)\int f\dif\mu.
  \end{equation}
  It follows from the proof of Greco's representation theorem (see   \cite[Thm.~13.2]{denneberg1994}) that we can actually assume $\mu$ to be defined as the   restriction of $\lowfnl_*$ to events:
  \begin{equation}\label{eq:definition-mu}
    \mu(A)
    =\lowfnl_*(A)
    =\sup\set{\lowfnl(f)}{f\leq I_A\text{ and }f\in\Domain\lowfnl}
  \end{equation}
  for all $A\subseteq\pspace$.  By Theorem~\ref{theo:natext-monotone-gambles-linlat}, $\mu$ is   also equal to the restriction to events of the natural extension $\lnex_\lowfnl=\lowfnl_*$   of $\lowfnl$.  Let us consider $A\subseteq B\subseteq\pspace$, and show that   $\lnex_\lowfnl(I_A+I_B)=\lnex_\lowfnl(I_A)+\lnex_\lowfnl(I_B)=\mu(A)+\mu(B)$.  Since the   exactness of $\lnex_\lowfnl$ implies that it is super-additive, we only need to prove that   $\lnex_\lowfnl(I_A+I_B)\leq\mu(A)+\mu(B)$. Given $\epsilon>0$, we deduce from   Eq.~\eqref{eq:inner-extension} that there is some $f$ in $\Domain\lowfnl$ such that $f\leq   I_A+I_B$ and $\lnex_\lowfnl(I_A+I_B)\leq\lowfnl(f)+\epsilon$.  Note that we may assume   without loss of generality that $f$ is non-negative [because $f\vee0$ belongs to   $\Domain\lowfnl$ and also satisfies the same inequality].  Let us define $g_1=f\wedge1$ and   $g_2=f-f\wedge1$.  These gambles belong to the linear lattice $\Domain\lowfnl$.  Moreover,   $g_1+g_2=f$. Let us show that $g_1\leq I_B$ and $g_2\leq I_A$:
  \begin{itemize}
  \item[--] Given $\omega\notin B$, we have $0\leq f(\omega)\leq(I_A+I_B)(\omega)=0$ whence     $g_1(\omega)=g_2(\omega)=0$.
  \item[--] Given $\omega \in A$, there are two possibilities: if $f(\omega)\leq 1$, then     $g_2(\omega)=0$ and $g_1(\omega)=f(\omega) \leq 1$. If on the other hand $f(\omega)> 1$,     then $g_1(\omega)=1$ and $g_2(\omega)=f(\omega)-1\leq 2-1=1$.
  \item[--] Given $\omega\in B\setminus A$, we have $f(\omega)\leq1$, whence     $g_1(\omega)=f(\omega)\leq1$ and $g_2(\omega)=0$.
  \end{itemize}
  Moreover, $g_1$ and $g_2$ are comonotone: consider any $\omega_1$ and $\omega_2$ in   $\pspace$, and assume that $g_2(\omega_1)<g_2(\omega_2)$. Then $g_2(\omega_2)>0$ and   consequently $\omega_2\in A$ and $f(\omega_2)>1$. This implies in turn that indeed   $g_1(\omega_2)=1\geq g_1(\omega_1)$.  Hence, since $\lowfnl$ is assumed to be comonotone   additive,
  \begin{equation*}
    \lnex_\lowfnl(I_A+I_B)\leq\lowfnl(f)+\epsilon
    =\lowfnl(g_1+g_2)+\epsilon
    =\lowfnl(g_1)+\lowfnl(g_2)+\epsilon \leq\lnex_\lowfnl(A)+\lnex_\lowfnl(B)+\epsilon,
  \end{equation*}
  and since this holds for all $\epsilon>0$ we deduce that indeed   $\lnex_\lowfnl(I_A+I_B)\leq\lnex_\lowfnl(A)+\lnex_\lowfnl(B)=\mu(A)+\mu(B)$.

  Now consider two arbitrary subsets $C$ and $D$ of $\pspace$. Then $C\cap D\subseteq C\cup   D$, and consequently
  \begin{multline*}
    \mu(C\cup D)+\mu(C\cap D) =\lnex_\lowfnl(I_{C\cup D}+I_{C \cap D}) =\lnex_\lowfnl(I_C+I_D)     \\ \geq \lnex_\lowfnl(I_C)+\lnex_\lowfnl(I_D) =\mu(C)+\mu(D),
  \end{multline*}
  taking into account that $\lnex_\lowfnl$ is super-additive (because it is exact). We   conclude that $\mu$ is $2$-monotone on $\wp(\pspace)$. From   Proposition~\ref{prop:charac-coherence}, we conclude that $\mu$ is an exact set function on   $\wp(\pspace)$, so by Theorem~\ref{theo:natext-choquet-lattice}, its natural extension is   the Choquet functional associated with $\mu$, and is therefore equal to $\lowfnl$ on   $\Domain\lowfnl$, by Eq.~\eqref{eq:equality-everywhere}. If we now apply   Theorem~\ref{theo:natext-monotone-gambles}, we see that the exact functional $\lowfnl$,   which has been shown to satisfy $\lowfnl(f)=(C)\int f\dif\mu$ for all $f$ in   $\Domain\lowfnl$, is also $2$-monotone.

  We now prove the converse implication. Assume that $\lowfnl$ is $2$-monotone.  Then,   applying Theorems~\ref{theo:inner-extension-monotone-gambles}   and~\ref{theo:natext-monotone-gambles-linlat}, its natural extension   $\lnex_\lowfnl=\lowfnl_*$ to all gambles is also $2$-monotone, and consequently so is its   restriction $\mu$ to events. Moreover, $\gambles_{\wp(\pspace)}=\gambles$, because any   gamble is the uniform limit of some sequence of simple gambles.  If we now apply   Theorem~\ref{theo:monotone-representation}, we see that $\lnex_\lowfnl(f)=(C)\int f\dif\mu$   for all $f$ in $\gambles$.  Consequently, $\lnex_\lowfnl$ is comonotone additive, because   the Choquet functional associated with a monotone set function is (see   \cite[Prop.~5.1]{denneberg1994}), and so is therefore $\lowfnl$.
\end{proof}

Hence, the natural extension of an $n$-monotone ($n\ge 2$) exact functional defined on a linear lattice of gambles that contains the constant gambles is always comonotone additive. Indeed, this natural extension is the Choquet functional associated to its restriction to events.

\begin{corollary}\label{cor:n-monolowprev:lnex:choqint}
  Let $n\in\SetNinfty$, $n\ge 2$, and let\/ $\lowfnl$ be an $n$-monotone exact functional   defined on a linear lattice that contains all constant gambles.  Then $\lnex_\lowfnl$ is   $n$-monotone, is comonotone additive, and is equal to the Choquet integral with respect to   $\lowfnl_*$ restricted to events.
\end{corollary}

Moreover, such an exact functional is generally not uniquely determined by its restriction to events, but it is uniquely determined by the values that its natural extension $\lnex_\lowfnl=\lowfnl_*$ assumes on events. Of course, this natural extension also depends in general on the values that $\lowfnl$ assumes on gambles, as is evident from Eq.~\eqref{eq:definition-mu}. On the other hand, we also deduce from the theorem that the procedure of natural extension preserves comonotone additivity from (indicators of) events to gambles.

As a nice side result, we deduce that an $n$-monotone ($n\ge 2$) exact set function $\lowfnl$ on $\wp(\pspace)$, which usually has many exact extensions to $\gambles$, has actually \emph{only one $2$-monotone} exact extension to $\gambles$. This unique $2$-monotone exact extension coincides with the natural extension of $\lowfnl$.

\begin{corollary}
  Let $n\in\SetNinfty$, $n\ge 2$. An $n$-monotone exact set function defined on all events has   a \emph{unique} $2$-monotone (or equivalently, comonotone additive) exact extension to all   gambles, that is furthermore automatically also $n$-monotone, namely its natural extension.
\end{corollary}
\begin{proof}
  Let $\lowfnl$ be an $n$-monotone exact set function defined on all events.  By   Theorem~\ref{theo:natext-monotone-gambles}, its natural extension $\lnex_\lowfnl$ to   $\gambles$ is an $n$-monotone, and hence, $2$-monotone exact extension of $\lowfnl$. The   proof is complete if we can show that $\lnex_\lowfnl$ is the only $2$-monotone exact   extension of $\lowfnl$.

  So, let $\alowfnl$ be any $2$-monotone exact extension of $\lowfnl$. We must show that   $\alowfnl=\lnex_\lowfnl$. Let $f$ be any gamble on $\pspace$. Then
  \begin{equation*}
    \alowfnl(f)
    =(C)\int f\dif\alowfnl
    =(C)\int f\dif\lowfnl
    =\lnex_\lowfnl(f),
  \end{equation*}
  where the first equality follows from Corollary~\ref{cor:n-monolowprev:lnex:choqint}, the   second from the equality of $\alowfnl$ and $\lowfnl$ on events, and the third by applying   Theorem~\ref{theo:natext-choquet-lattice}.  This establishes uniqueness.
\end{proof}

We summarise some of the comments and results in this section in Figure~\ref{fig:relationships}.

\begin{figure}[htbp]
  \begin{center}
    \begin{equation*}
      \begin{CD}
        @. \text{$\lowfnl$ $n$-monotone} @>>> \text{$\alowfnl$ $n$-monotone}\\
        @. @VVV @VVV\\
        \text{$\lowfnl$ comonotone additive} @= \text{$\lowfnl$ $2$-monotone}
        @>>> \text{$\alowfnl$ $2$-monotone}\\
        @. @| @|\\
        @. \lowfnl=\lnex_\alowfnl=(C)\int\cdot\dif\alowfnl @>>>         \lnex_\alowfnl=(C)\int\cdot\dif\alowfnl
      \end{CD}
    \end{equation*}
    \caption{Relationships between the properties of an exact functional $\lowfnl$ on       $\gambles$ and its restriction $\alowfnl$ to events; implications are depicted using       arrows, equivalences using double lines.}
    \label{fig:relationships}
  \end{center}
\end{figure}

Next, we relate comonotone additivity, or equivalently, $2$-monotonicity, of exact functionals to properties of their sets of dominating linear exact functionals.

\begin{proposition}\label{prop:2-monotone-linear-gambles}
  Let\/ $\lowfnl$ be an exact functional on a linear lattice of gambles.
  \begin{enumerate}[(a)]
  \item If\/ $\lowfnl$ is comonotone additive on its domain, then for all comonotone $f$ and     $g$ in $\Domain\lpr$, there is some $\alinfnl$ in $\solp(\lowfnl)$ such that\/     $\alinfnl(f)=\lowfnl(f)$ and\/ $\alinfnl(g)=\lowfnl(g)$.
  \item Assume in addition that\/ $\Domain\lowfnl$ contains all constant gambles.  Then     $\lowfnl$ is comonotone additive (or equivalently $2$-monotone) on its domain if and only     if for all comonotone $f$ and $g$ in $\Domain\lowfnl$, there is some $\alinfnl$ in     $\solp(\lowfnl)$ such that\/ $\alinfnl(f)=\lowfnl(f)$ and\/ $\alinfnl(g)=\lowfnl(g)$.
  \end{enumerate}
\end{proposition}

\begin{proof}
  To prove the first statement, assume that $\lowfnl$ is comonotone additive on its domain,   and consider $f$ and $g$ in $\Domain\lowfnl$ that are comonotone.  Then $f+g$ also belongs   to $\Domain\lowfnl$, so we know that $\lowfnl(f+g)=\lowfnl(f)+\lowfnl(g)$. On the other   hand, since $\lowfnl$ is exact, there is some $\alinfnl$ in $\solp(\lowfnl)$ such that   $\lowfnl(f+g)=\alinfnl(f+g)=\alinfnl(f)+\alinfnl(g)$. So   $\alinfnl(f)+\alinfnl(g)=\lowfnl(f)+\lowfnl(g)$ and since we know that   $\lowfnl(f)\leq\alinfnl(f)$ and $\lowfnl(g)\leq\alinfnl(g)$, this implies that   $\lowfnl(f)=\alinfnl(f)$ and $\lowfnl(g)=\alinfnl(g)$.

  The `only if' part of the second statement is an immediate consequence of the first. To   prove the `if' part, consider arbitrary comonotone $f$ and $g$ in $\Domain\lowfnl$. Then it   is easy to see that $f\join g$ and $f\meet g$ are comonotone as well, and belong to   $\Domain\lowfnl$, so by assumption there is a $\alinfnl$ in $\solp(\lowfnl)$ such that   $\lowfnl(f\meet g)=\alinfnl(f\meet g)$ and $\lowfnl(f\join g)=\alinfnl(f\join g)$.  Then,   using Theorem~\ref{theo:linprevs-completely-everything},
  \begin{equation*}
    \lowfnl(f\join g)+\lowfnl(f\meet g)=\alinfnl(f\join g)+\alinfnl(f\meet g)
    =\alinfnl(f)+\alinfnl(g)\geq\lowfnl(f)+\lowfnl(g).
  \end{equation*}
  This tells us that $\lowfnl$ is $2$-monotone, and by Theorem~\ref{theo:equivalence} also   comonotone additive.
\end{proof}

As a corollary, we deduce the following, apparently first proven by Walley \cite[Cors.~6.4 and~6.5, p.~57]{walley1981} for coherent lower previsions.

\begin{corollary}\label{cor:2-monotone-linear-events}
  Let\/ $\lowfnl$ be an exact set function on a lattice of events.  Then $\lowfnl$ is   $2$-monotone if and only if for all $A$ and $B$ in $\Domain\lowfnl$ such that $A\subseteq   B$, there is some $\alinfnl$ in $\solp(\lowfnl)$ such that\/ $\alinfnl(A)=\lowfnl(A)$ and\/   $\alinfnl(B)=\lowfnl(B)$.
\end{corollary}

\begin{proof}
  We just show that the direct implication is a consequence of the previous results; the   converse one follows easily by applying the condition to $A\cap B\subseteq A\cup B$, for $A$   and $B$ in $\Domain\lowfnl$.

  Let $\lowfnl$ be a $2$-monotone exact functional defined on a lattice of events.  From   Theorem~\ref{theo:natext-monotone-gambles}, the natural extension $\lnex_\lowfnl$ of   $\lowfnl$ to all gambles is $2$-monotone and exact. Hence, given $A\subseteq   B\in\Domain\lowfnl$, since $I_A$ and $I_B$ are comonotone,   Proposition~\ref{prop:2-monotone-linear-gambles} implies the existence of a $\alinfnl$ in   $\solp(\lnex_\lowfnl)=\solp(\lowfnl)$ such that   $\alowfnl(A)=\lnex_\lowfnl(A)=\lowfnl_*(A)=\lowfnl(A)$ and   $\alowfnl(B)=\lnex_\lowfnl(B)=\lowfnl_*(B)=\lowfnl(B)$.
\end{proof}

\section{Conclusions}\label{sec:conclusions}

We see from the results in this paper that there is no real reason to restrict the notion of $n$-monotonicity to set functions (or lower probabilities). In fact, it turns out that it is fairly easy, and completely within the spirit of Choquet's original definition, to define and study this property for functionals (or lower previsions). And in fact, we have shown above that doing this does not lead to just another generalisation of something that existed before, but that it leads to genuinely new insights. Our results also show that the procedure of natural extension is of particular interest for $n$-monotone lower previsions; not only does it provide the behaviourally most conservative (i.e., point-wise smallest) extension to all gambles, but it is also the only extension to be $n$-monotone: hence, any other extension is implying behavioural dispositions that are not implied by coherence (alone), \emph{and at the   same time} it does not satisfy $2$-monotonicity.

We deduce from our results that, under exactness, $2$-monotonicity of a lower prevision is actually equivalent to comonotone additivity, and therefore to being representable as a Choquet functional (see Theorem~\ref{theo:equivalence} for a precise formulation). In particular, this means that all the results we have established in this paper for $2$-monotone exact functionals are valid for comonotone additive functionals.

Finally, we would like to mention that we have shown elsewhere (\cite{cooman2006}) that most (if not all) of the lower integrals defined in the literature are actually completely monotone, and are therefore representable as a Choquet functional. Indeed, we also show in that paper that we can use most of the lower integrals in the literature to calculate the natural extension of bounded charges, and of some finitely additive set functions.

\section*{Acknowledgements}

This paper has been partially supported by research grant G.0139.01 of the Flemish Fund for Scientific Research (FWO), the Belgian American Educational Foundation and by the projects MTM2004-01269, TSI2004-06801-C04-01.


\begin{thebibliography}{10}

\bibitem{aigner1977} M.~Aigner.  \newblock {\em Combinatorial Theory}.  \newblock Classics in   Mathematics. Springer-Verlag, Berlin, 1977.

\bibitem{artzner1999} Ph. Artzner, F.~Delbaen, J.-M. Eber, and D.~Heath.  \newblock Coherent   measures of risk.  \newblock {\em Mathematical Finance}, 9:203--228, 1999.

\bibitem{bhaskara1983} K.~P.~S. Bhaskara~Rao and M.~Bhaskara~Rao.  \newblock {\em Theory of     Charges}.  \newblock Academic Press, London, 1983.

\bibitem{choquet1953} G.~Choquet.  \newblock Theory of capacities.  \newblock {\em Annales de     l'Institut {F}ourier}, 5:131--295, 1953--1954.

\bibitem{cooman2001} G.~{d}e Cooman.  \newblock Integration and conditioning in numerical   possibility theory.  \newblock {\em Annals of Mathematics and Artificial Intelligence},   32:87--123, 2001.

\bibitem{cooman1998} G.~{d}e Cooman and D.~Aeyels.  \newblock Supremum preserving upper   probabilities.  \newblock {\em Information Sciences}, 118:173--212, 1999.

\bibitem{cooman1998a} G.~{d}e Cooman and D.~Aeyels.  \newblock A random set description of a   possibility measure and its natural extension.  \newblock {\em IEEE Transactions on Systems,     Man and Cybernetics---Part A: Systems and Humans}, 30:124--130, 2000.

\bibitem{cooman2006} G.~{d}e Cooman, M.~Troffaes, and E.~Miranda.  \newblock A unifying   approach to integration for bounded positive charges.  \newblock Submitted for publication,   2006.

\bibitem{delbaen2002} F.~Delbaen.  \newblock Coherent risk measures on general probability   spaces.  \newblock In K.~Sandmann and P.~J. Sch\"onbucher, editors, {\em Advances in Finance     and Stochastics}, pages 1-- 37. Delbaen, F, Berlin, 2002.

\bibitem{denneberg1994} D.~Denneberg.  \newblock {\em Non-Additive Measure and Integral}.   \newblock Kluwer Academic, Dordrecht, 1994.

\bibitem{dubois1988} D.~Dubois and H.~Prade.  \newblock {\em Possibility Theory}.  \newblock   Plenum Press, New York, 1988.

\bibitem{2003:ferson::pboxes} Scott Ferson, Vladik Kreinovich, Lev Ginzburg, Davis~S. Myers,   and Kari Sentz.  \newblock Constructing probability boxes and {D}empster-{S}hafer   structures.  \newblock Technical Report SAND2002--4015, Sandia National Laboratories,   January 2003.

\bibitem{levi1980a} I.~Levi.  \newblock {\em The Enterprise of Knowledge}.  \newblock MIT   Press, London, 1980.

\bibitem{maass2002} S.~Maa\ss.  \newblock Exact functionals and their core.  \newblock {\em     Statistical Papers}, 43:75--93, 2002.

\bibitem{maass2003} S.~Maa\ss.  \newblock {\em Exact functionals, functionals preserving     linear inequalities, L\'evy's metric}.  \newblock PhD thesis, University of Bremen, 2003.

\bibitem{nguyen1997} H.~T. Nguyen, N.~T. Nguyen, and T.~Wang.  \newblock On capacity   functionals in interval probabilities.  \newblock {\em International Journal of Uncertainty,     Fuzziness and Knowledge-Based Systems}, 5:359--377, 1997.

\bibitem{schechter1997} E.~Schechter.  \newblock {\em Handbook of Analysis and Its     Foundations}.  \newblock Academic Press, San Diego, CA, 1997.

\bibitem{schmeidler1972} D.~Schmeidler.  \newblock Cores of exact games.  \newblock {\em     Journal of Mathematical Analysis and Applications}, 40:214--225, 1972.

\bibitem{walley1981} P.~Walley.  \newblock Coherent lower (and upper) probabilities.   \newblock Technical report, University of Warwick, Coventry, 1981.  \newblock Statistics   Research Report 22.

\bibitem{walley1991} P.~Walley.  \newblock {\em Statistical Reasoning with Imprecise     Probabilities}.  \newblock Chapman and Hall, London, 1991.

\end{thebibliography}
\end{document}